\documentclass[11pt]{article}
\usepackage{authblk}
\usepackage{amssymb}
\usepackage{amsthm}
\usepackage{graphicx}
\usepackage{enumerate}
\usepackage{mathtools}
\usepackage{fullpage}
\usepackage{xcolor}
\usepackage{comment}
\usepackage{tikz}
\usepackage{bbm}
\usetikzlibrary{patterns}
\usetikzlibrary{decorations.pathmorphing}
\usetikzlibrary{decorations.pathreplacing,decorations.markings}
\usetikzlibrary{shapes}

\tikzset{every state/.style={minimum size=0pt}}

\tikzset{
	mark position/.style args={#1(#2)}{
		postaction={
			decorate,
			decoration={
				markings,
				mark=at position #1 with \coordinate (#2);
			}
		}
	}
}

\usetikzlibrary{3d}
\usepackage{xcolor}
\usepackage{xifthen}
\usepackage{blindtext}
\usepackage{geometry}
\newtheorem{thm}{Theorem}[section]
\newtheorem{lm}[thm]{Lemma}

\newtheorem{prop}[thm]{Proposition}

\newtheorem{rmk}[thm]{Remark}

\numberwithin{equation}{section}

\newcommand{\bigO}{\mathcal{O}}

\newcommand*\dif{\mathop{}\!\mathrm{d}}
\newcommand{\ddbar}[2]{\frac{{\mathrm d}#1}{2\pi {\mathrm i}#2}}

\newcommand{\FGUE}{\mathrm{F}_{\mathrm{GUE}}}

\newcommand{\ii}{\mathrm{i}}
\newcommand{\inn}{\mathrm{in}}

\newcommand{\LL}{\mathrm{L}}
\newcommand{\EE}{\mathbb{E}}
\newcommand{\out}{\mathrm{out}}
\newcommand{\prob}{\mathbb{P}}
\newcommand{\realR}{\mathbb{R}}

\newcommand{\RR}{\mathrm{R}}

\newcommand{\rz}{\mathrm{z}}

\newcommand{\rh}{\mathrm{h}}
\newcommand{\rt}{\mathrm{T}}

\newcommand{\Ai}{\mathrm{Ai}}
\newcommand{\HH}{\mathcal{H}}
\newcommand{\Addresses}{{% additional braces for segregating \footnotesize
  \bigskip
  \footnotesize

  Ron Nissim, \textsc{Department of Mathematics, Massachusetts Institute of Technology, Cambridge, Massachusetts 02139, USA.}\par\nopagebreak
  \textit{E-mail address}: \texttt{rnissim@mit.edu}

  \medskip

  Ruixuan Zhang, \textsc{Department of Mathematics, University of Kansas, Lawrence, Kansas 66045, USA.}\par\nopagebreak
  \textit{E-mail address}: \texttt{rayzhang@ku.edu}
}}

\begin{document}

\title{Edgeworth-type expansion for the one-point distribution of the KPZ fixed point with a large height at a prior location
%The limiting conditional one-point distribution of the KPZ fixed point on a plateau
}

\date{}
\author{Ron Nissim, Ruixuan Zhang}
\maketitle

\begin{abstract}
 We consider the Kardar-Parisi-Zhang (KPZ) fixed point $\mathrm{H}(x,\tau)$ with the narrow-wedge initial condition and investigate the distribution of $\mathrm{H}(x,\tau)$ conditioned on a large height at an earlier space-time point $\mathrm{H}(x',\tau')$. As $\mathrm{H}(x',\tau')$ tends to infinity, we prove that the conditional one-point distribution of $\mathrm{H}(x,\tau)$ in the regime $\tau>\tau'$ converges to the Gaussian Unitary Ensemble (GUE) Tracy-Widom distribution and that the next two lower-order error terms can be expressed as derivatives of the Tracy-Widom distribution. The lowe order expansion here is analogue to the Edgeworth expansion in the central limit theorem. These KPZ-type limiting behaviors are different from the Gaussian-type ones obtained in \cite{Liu-Wang22} where they study the finite-dimensional distribution of $\mathrm{H}(x,\tau)$ conditioned on a large height at a later space-time point $\mathrm{H}(x',\tau')$. They show, with the narrow-wedge initial condition, that the conditional random field $\mathrm{H}(x,\tau)$ in the regime $\tau<\tau'$ converges to the minimum of two independent Brownian bridges modified by linear drifts as $\mathrm{H}(x',\tau')$ goes to infinity. The two results stated above provide the phase diagram of the asymptotic behaviors of a conditional law of KPZ fixed point in the regimes $\tau>\tau'$ and $\tau<\tau'$ when $\mathrm{H}(x',\tau')$ goes to infinity.
\end{abstract}

\section{Introduction}
\subsection{Background on the KPZ universality class and the KPZ fixed point.}

The Kardar-Parisi-Zhang (KPZ) equation is a nonlinear stochastic partial differential equation (PDE) given by 
\begin{equation}
    \partial_t \HH=\lambda\partial^2_x \HH +\nu\left(\partial_x \HH \right)^2+\xi,
\end{equation}
 where $\HH:(0,\infty)\times\realR\to \realR$ and $\xi:(0,\infty)\times \realR \to \realR $ represents spacetime white noise. It was first introduced by Kardar, Parisi, and Zhang in their seminal work \cite{KPZ86} to describe the evolution of the interface in a wide range of random growth models.
 Numerous papers in physics literature \cite{FNS77,vBKS85} had previously predicted that the height function $\HH(t,x)$, when appropriately scaled, should converge to a universal random field that is independent of the specific model. Over the past two decades, an extensive class of (1 + 1) spacetime dimension models, collectively referred to as the KPZ universality class, has been exactly solved. Remarkably, all of these models have been shown to share the same scaling limit. Some notable examples include the longest increasing subsequence \cite{BDJ}, directed last passage percolation \cite{Jo00}, polynuclear growth models \cite{Jo03}, the asymmetric simple exclusion process (ASEP) \cite{TW08,TW09,Virag20,QS22}, directed random polymers \cite{ACQ,S12} and the KPZ equation itself \cite{QS22,Virag20}.

The limit space-time field $\mathrm{H}(x,\tau)$, where $x\in \realR, t\ge 0$, of the KPZ universality class, is called the KPZ fixed point, which depends on the initial condition $\mathrm{H}(x,0)=h_0(x)$ for a function $h_0$ in the space of upper-semi-continuous functions.  It was first constructed by Matetski, Quastel and Remenik in \cite{MQR} recently, as a Markov process with explicit transition probability by analyzing the totally asymmetric simple exclusion process (TASEP). However, the derivation of explicit formulas for the multi-point distribution of $\mathrm{H}(x,\tau)$ is a challenging task due to the complexity of the transition probability formula. Pioneering efforts in this field have unfolded over the past two decades. Notably, \cite{BDJ,Jo00} demonstrated that the marginal one-point distribution of $\mathrm{H}(x,\tau)$ follows the Tracy-Widom distribution. Meanwhile, the multi-point distribution along the spatial direction, $\left(\mathrm{H}(x_1,\tau),\cdots,\mathrm{H}(x_n,\tau)\right)$ for a fixed $\tau$, is characterized by the finite-dimensional distribution of the Airy process and its analogues, as explored in works by \cite{PS02,Jo03,IS04,BFP07,BFPS,BFS08,MQR}. Recently, significant progress has been made in obtaining general explicit formulas for joint distributions, even when considering different temporal points $\left(\mathrm{H}(x_1,\tau_1),\cdots,\mathrm{H}(x_n,\tau_n)\right)$. These breakthroughs were achieved by Johansson and Raham \cite{JR21} as well as independently by Liu \cite{Liu19}.

The multi-point distribution formula enables us to delve into the finer details of the KPZ fixed point. For an inhomogeneous Markov process, we naturally turn our attention to the probabilistic properties of the conditional field $\{\mathrm{H}(x,\tau)\mid \mathrm{H}(x',\tau')=\rh'\}$ for given $x'\in \realR,\tau'>0$. The primary goal is to provide insights into its asymptotic behavior when the height $\rh'$ grows to infinity. To achieve this, we concentrate on the narrow wedge (step) initial condition $\mathrm{H}(x,0)=-\infty \mathbf{1}_{x\ne 0}$.

We discover that when conditioned on this rare event, the behavior of the KPZ fixed point $\mathrm{H}(x,\tau)$ exhibits distinct characteristics for two scenarios: $\tau<\tau'$ and $\tau>\tau'$. In a recent study by \cite{Liu-Wang22}, it was demonstrated that when $\tau<\tau'$, they have, as $\rh'\to \infty$, 

\begin{equation}
\label{eq:KPZ_Brownian}
    \begin{split}
         &\mathrm{Law}\left(\left\{\left.\frac{\mathrm{H}\left(\frac{x\tau'^{3/4}}{\sqrt{2}\rh'^{1/4}}+tx',t\tau'\right)-t\mathrm{H}(x',\tau')}{\sqrt{2}\left(\tau'\rh'\right)^{1/4}}\right\}_{x\in\realR,t\in(0,1)}\right|\mathrm{H}(x',\tau')=\rh'\right)\\
         &\xrightarrow{f.d.d.} \mathrm{Law}\left(\left\{\min\left\{\mathbb{B}_1(t)+x,\mathbb{B}_2(t)-x\right\}\right\}_{x\in\realR,t\in(0,1)}\right)
    \end{split}
\end{equation}
 where $\xrightarrow{f.d.d.}$ denotes convergence of finite-dimensional distributions. Here $\left\{\mathbb{B}_1(t)\right\}_{t\in (0,1)}$ and $\left\{\mathbb{B}_2(t)\right\}_{t\in (0,1)}$ are two independent Brownian bridges on $[0,1]$.

In our current paper, we primarily focus on the second case, where $\tau>\tau'$. In this context, we anticipate that the KPZ fixed point, after this high point, asymptotically behaves in a manner consistent with the conventional KPZ fixed point, as detailed in Remark~\ref{rmk:01}. We will validate this by examining the one-point distribution. Specifically, we will show that the  $\{\mathrm{H}(x,\tau)-\mathrm{H}(x',\tau')\}$, conditioned on a large value of $\mathrm{H}(x',\tau')$, follows the GUE Tracy-Widom distribution, complemented by smaller-order terms associated with the derivatives of the GUE Tracy-Widom distribution. This result is elucidated in Theorem~\ref{thm:main}. The case where $\tau=\tau'$ is not addressed in \cite{Liu-Wang22} or our paper for distinct reasons. Notably, a phase transition occurs at $\tau=\tau'$ regarding the asymptotic behavior of $\{\mathrm{H}(x,\tau)\mid \mathrm{H}(x',\tau')=\rh'\}$. In an upcoming work \cite{LiuZhang24}, the authors have established a limit theorem for the conditional field $\{\mathrm{H}(x,\tau)\mid \mathrm{H}(0,1)=L\}$ in the scenario where $\tau=\tau_L\to 1$ following an appropriate scaling as $L\to\infty.$ The transition from Gaussian Universality to KPZ Universality is frequently observed in the study of the KPZ universal class, as reviewed in \cite{IC16}.

In Theorem \ref{thm:main}, we present the explicit formulations of the second and third-order expansions, which are conceptually similar to the classical Edgeworth expansion. These expansions contribute asymptotic correction terms to the Central Limit Theorem up to an order that depends on the number of moments available. Roughly speaking, if $\{X_i\}_{i=1}^n$ is a sequence of i.i.d. random variables with $\EE |X_i|^{j+2}<\infty$ and $\limsup_{|t|\to \infty} \left\|\EE e^{\ii tX_i}\right\|<\infty$, then Edgeworth expansion of the cumulative distribution function of $S_n:=\frac{\sqrt{n}\left(\sum_{i=1}^n X_i-\mu\right)}{\sigma}$ can be written as 
\begin{equation*}
    \prob\left(S_n\leq x\right)=\Phi(x)+n^{-1/2}p_1(x)\phi(x)+n^{-1}p_2(x)\phi(x)+\dots+n^{-j/2}p_j(x)\phi(x)+o(n^{-j/2})
\end{equation*}
uniformly in all $x$, where $\phi(x)=(2\pi)^{-1/2}e^{-x^2/2}$ is the standard normal density function and $\Phi(x)=\int_{-\infty}^x (2\pi)^{-1/2}e^{-t^2/2}dt$ is the standard normal distribution function. The functions $p_j(x)$ are Hermite polynomials with coefficients depending on the cumulants of $X_i$. For a detailed discussion on this topic, we refer to the work of \cite{wallace1958asymptotic,hall2013bootstrap}. Notably, the explicit expression of error terms in the context of the KPZ universality class is a rarity in the literature, while the coefficients are expected to depend on different models. It will be of great significance to explore this phenomenon in the context of specific integrable models, such as TASEP, last passage percolation and directed random polymers.

Moreover, addressing the scenario where $\rh' \to -\infty$ poses a significant challenge. The asymptotic analysis of this problem would be considerably more complex, and we have opted to defer this investigation to future research.

\subsection{Main results.}
Let
\begin{equation}
    \mathrm{F}(\rh; x,\tau):= \prob\left(\mathrm{H}(x,\tau) \le \rh\right)=\FGUE\left(\frac{1}{\tau^{1/3}}\rh+\frac{1}{\tau^{4/3}}x^2\right), \label{FGUE}
\end{equation}
where $\FGUE$ is the cumulative distribution function of the GUE Tracy-Widom distribution.  It is the distribution of the normalized largest eigenvalue of a random Hermitian matrix \cite{TW96}. Its cumulative distribution function can be given as an integral \cite{TW96},
\begin{equation}
    \FGUE(x)=\exp{\left(-\int_x^{\infty}(s-x)q^2(s)ds\right)}.
\end{equation}
Here, the function $q(x)$ is the solution to the Painlev\'e II equation
\begin{equation}
    q''=2q^3+xq
\end{equation}
with the boundary condition $q(x)\sim \Ai(x)$, $x\to \infty,$ where $\Ai(x)$ is the Airy function satisfying $\Ai''(x)=x\Ai(x).$ We used the notation $f(x)\sim g(x)$ as $x\to c$  if $f(x)/g(x) \to 1$ as $x\to c$. It is well-known that $\FGUE(x)$ has asymptotics
\begin{equation}
    1-\FGUE(x)\sim \frac{1}{16\pi x^{3/2}}\exp{\left(-\frac{4}{3}x^{3/2}\right)} \label{TWasymp1}
\end{equation}
as $x\to \infty,$ and its density function
\begin{equation}
    \mathrm{F}'_{\mathrm{GUE}}(x)\sim\frac{1}{8\pi x}\exp{\left(-\frac{4}{3}x^{3/2}\right)} \label{TWasymp2}
\end{equation} as $x\to \infty.$ 
See \cite{BBD08} for more details.

We will now present our main results. 
\begin{thm}
    \label{thm:main}
      Suppose $x,x',\rh\in \mathbb{R}$ and $\tau,\tau'>0$ are fixed constants. Then
     \begin{equation}
     \label{eq:main_result}
         \begin{split}
             	&\prob\left(\mathrm{H}(x+x',\tau+\tau')-\mathrm{H}(x',\tau')\le \rh\mid \mathrm{H}(x',\tau')=\rh'\right)\\
             	&=\mathrm{F}(\rh;x,\tau)-\frac{3}{2}\sqrt{\frac{\tau'}{\rh'}}\mathrm{F}'(\rh;x,\tau)+\frac{\tau'}{\rh'}\mathrm{F}''(\rh;x,\tau)+\bigO((\rh')^{-3/2})
         \end{split}
     \end{equation}
     as $\rh'\to \infty$, where $\mathrm{F}'(\rh;x,\tau)$ and $\mathrm{F}''(\rh;x,\tau)$ stand for the first and second order derivatives of $\mathrm{F}(\rh;x,\tau)$. Here the conditional probability $\prob(A \mid \mathrm{H}(x',\tau')=\rh')$ is understood as $\lim_{\epsilon\to 0}\prob(A \mid |\mathrm{H}(x',\tau')-\rh'|<\epsilon)$ for any measurable event $A$.
\end{thm}

\begin{rmk}
\label{rmk:01}
The leading term of equation \eqref{eq:main_result} being $\mathrm{F}(\rh;x,\tau)$ is not unexpected. We can interpret the KPZ fixed point $\mathrm{H}(x,\tau)$ in the language of the directed landscape $\mathcal{L}(y,s;x,t)$ \cite{DOV}, which was introduced as the random field arising in the limit of Brownian last passage
percolation \cite{DOV} and was proved or expected to be the universal limit of the KPZ universality class \cite{Dau21}. For the step initial condition, the KPZ fixed point and the directed landscape are related by the following equation  \cite{NQR20},
\begin{equation}
    \mathrm{H}(x,\tau)=\mathcal{L}(0,0;x,\tau).
\end{equation}

One basic property of the directed landscape $\mathcal{L}$ is the following composition law
\begin{equation}
    \mathcal{L}(y,s;x,t)=\max_{z\in\realR} \{\mathcal{L}(y,s;z,r)+\mathcal{L}(z,r;x,t)\}
\end{equation}
for any $r$ satisfying $s<r<t$. Especially,
\begin{equation}
\label{eq:identity_landscape}
    \mathcal{L}(0,0;x+x',\tau+\tau')=\max_{y\in\realR}\{\mathcal{L}(0,0;y,\tau')+\mathcal{L}(y,\tau';x+x',\tau+\tau')\}.
\end{equation}
Note that if we condition on the event $\mathcal{L}(0,0;x',\tau')=\mathrm{H}(x',\tau')$ becomes very large, it is reasonable to expect that the maximum on the right hand side of~\eqref{eq:identity_landscape} is achieved near the point $y=x'$. Thus we expect the equation~\eqref{eq:identity_landscape} becomes
\begin{equation}
    \mathrm{H}(x+x',\tau+\tau')\approx\mathrm{H}(x',\tau')+\mathcal{L}(x',\tau';x+x',\tau+\tau').
\end{equation}
Using the translation invariance property of the directed landscape \cite{DOV}, $\mathcal{L}(x',\tau';x+x',\tau+\tau')$ has the same distribution as $\mathcal{L}(0,0;x,\tau)$, which has $\mathrm{F}(\rh;x,
\tau)$ as its one-point distribution.
\end{rmk}
\begin{rmk}
The emergence of $\mathrm{F}'(\rh)$ and $\mathrm{F}''(\rh)$ as exact first and second-order error terms in the asymptotic analysis is indeed remarkable.  To the best of our knowledge, this result represents the most precise characterization of the conditional law of $\mathrm{H}(x+x',\tau+\tau')$ as the height of $\mathrm{H}(x',\tau')$ tends to infinity. This finding also extends to the asymptotic behavior of the two-point distribution of the KPZ fixed point, as detailed in Theorem \ref{twopointasy} and discussed further in Remark 1.6.

It may pique one's curiosity whether a discernible pattern exists for higher-order terms. In fact, the third-order term does involve $\mathrm{F}'''(\rh)$ alongside other complex terms. Regrettably, these additional terms cannot be expressed solely as derivatives of $\mathrm{F}(\rh)$ and are characterized by intricate complexity, lacking the simplicity of expression found in the earlier terms.
\end{rmk}

Furthermore, we are able to derive the asymptotic behavior of the two-point tail probability for the KPZ fixed point when it starts with a step initial condition. To arrive at this result, we employ a similar argument as that employed in Theorem \ref{thm:main}, focusing on \eqref{eq2} as defined in Proposition \ref{prop1}. Given the parallel nature of this derivation, we refrain from providing a detailed exposition here.

\begin{thm}\label{twopointasy}
Suppose $x,x',\rh\in \mathbb{R}$ and $\tau,\tau'>0$ are fixed constants. Then
		\begin{equation}
			\begin{split}
		    &\prob\left(\mathrm{H}(x',\tau')\ge \rh',\mathrm{H}(x+x',\tau+\tau')\ge \rh+\rh'\right)\\
			=&\prob\left(\mathrm{H}(x',\tau')\ge \rh'\right)\left(\prob\left(\mathrm{H}(x,\tau)\ge \rh\right)+2\sqrt{\frac{\tau'}{\rh'}}\mathrm{F}'(\rh;x,\tau)-2\frac{\tau'}{\rh'}\mathrm{F}''(\rh;x,\tau)+\bigO\left(\rh'^{-3/2}\right)\right)
			\end{split}
		\end{equation}
	as $\rh'\to\infty$. 
\end{thm}

\begin{rmk}
There is relevant work by \cite{GH22}, where they obtain several sharp results concerning the upper tail behavior of the KPZ equation. In particular, they established a precise tail estimate for the two-point distribution for all time values $t$ within the interval $[t_0,\infty)$ given by:
$$\prob\left(h(-\theta^{1/2},t)\ge a\theta,h(\theta^{1/2},t)\ge b\theta\right)=\exp{\left(-\frac{4}{3}\theta^{3/2}\left((1+a)^{3/2}+(1+b)^{3/2}\right)+\text{error}\right)}.$$ 
Here, $h(x,t):=t^{-1/3}(\HH (t,t^{2/3}x)+\frac{t}{12})$ is the solution of the KPZ equation $ \partial_t \HH=\frac{1}{4}\partial^2_x \HH +\frac{1}{4}\left(\partial_x \HH \right)^2+\xi.$ The error term comes with explicit bounds, as further elaborated in their work \cite[Theorem 7.1]{GH22}. The leading term they obtained matches the asymptotics of the right tail of the GUE Tracy–Widom distribution function, as shown in \eqref{TWasymp1}. Here, we provide the exact form of the leading term and lower-order terms, even for the upper tail of the two-time distribution.
\end{rmk}

\begin{rmk}
    In the physics community, \cite{de2017tail,de2018two} studies the joint probability distribution function of the height of the KPZ equation with wedge initial conditions, at two different times $t_1<t_2$ in the limit where both times are large and their ratio is fixed. Using the replica Bethe ansatz method, they obtain its exact tail when the height at the earlier time is large and positive. Their formula interpolates between two limits where the joint probability distribution function decouples: (i) into a product of two GUE Tracy–Widom (TW) distributions when $t_2/t_1\to\infty$, and (ii) into a product of a GUE-TW distribution and a Baik–Rains distribution when $(t_2-t_1)/t_1 \ll 1$. This result shows excellent agreement with experimental and numerical data \cite{de2017memory}, and also matches the rigorous expression for the joint distribution obtained in \cite{johansson2019two}. This decoupling phenomenon is similar to what is described in Theorem \eqref{twopointasy}, where, instead, we fix the time but take the height $\rh'\to \infty$, providing useful intuition for studying other models in the KPZ universality class.
\end{rmk}

\section*{Acknowledgements}
We would like to thank Jinho Baik, Zhipeng Liu and Yizao Wang for the comments and suggestions. R. Z. was partially supported by the NSF grant DMS-1953687 and Simons Collaboration Grant No. 637861. R.N. was supported by the NSF grant DMS-1954790. 

\section{Preliminaries}

We can write
\begin{equation}
    \begin{split}
        &\prob\left(\mathrm{H}(x+x',\tau+\tau')-\mathrm{H}(x',\tau')\le \rh\mid \mathrm{H}(x',\tau')=\rh'\right)\\
        &=1-\prob\left(\mathrm{H}(x+x',\tau+\tau')-\mathrm{H}(x',\tau')\ge \rh\mid \mathrm{H}(x',\tau')=\rh'\right).
    \end{split}
\end{equation}
Indeed, our focus narrows down to the investigation of $\prob\left(\mathrm{H}(x+x',\tau+\tau')-\mathrm{H}(x',\tau')\ge \rh\mid \mathrm{H}(x',\tau')=\rh'\right)$ as the primary concern. To reiterate, the conditional probability of interest can be examined through the following approach:
\begin{equation}
    \begin{split}
        &\prob\left(\mathrm{H}(x+x',\tau+\tau')-\mathrm{H}(x',\tau')\ge \rh\mid \mathrm{H}(x',\tau')=\rh'\right)\\
        &:=\lim_{\epsilon\to 0}\frac{\prob\left(\mathrm{H}(x+x',\tau+\tau')-\mathrm{H}(x',\tau')\ge \rh, \mathrm{H}(x',\tau')\in(\rh'-\epsilon,\rh'+\epsilon)\right)}{\prob(\mathrm{H}(x',\tau')\in(\rh'-\epsilon,\rh'+\epsilon))}\\
        &=\frac{\frac{\partial}{\partial\rh'}\prob\left(\mathrm{H}(\hat{x},\hat\tau)\ge \hat\rh,\mathrm{H}(x',\tau')\ge \rh'\right)\bigg|_{\hat{h}=h+h', \hat{x}=x+x', \hat{\tau}=\tau+\tau'}}{\frac{\mathrm{d}}{\mathrm{d}\rh'}\prob\left(\mathrm{H}(x',\tau')\ge \rh'\right)}.
    \end{split}\label{conditional}
\end{equation}

To advance our analysis, we initiate by referencing the formulas provided in \cite{Liu19} pertaining to the one-point and two-point distributions of $\mathrm{H}(x,\tau)$.

\subsection{Explicit formula for the joint distribution.}

We use the notations
\begin{equation}
    \begin{split}
        &\Delta(W):=\prod_{1 \leq i<j \leq n}\left(w_{j}-w_{i}\right), \quad \Delta\left(W ; W^{\prime}\right):=\prod_{i=1}^{n} \prod_{i^{\prime}=1}^{n^{\prime}}\left(w_{i}-w_{i^{\prime}}^{\prime}\right),\\ &f(W)=\prod_{i=1}^{n} f(w), \quad C(W;W'):=\frac{\Delta(W)\Delta(W')}{\Delta(W;W')}
    \end{split}
     \label{noation}
\end{equation}
for any two vectors $W=\left(w_{1}, \cdots, w_{n}\right) \in \mathbb{C}^{n}$ and $W^{\prime}=\left(w_{1}^{\prime}, \cdots, w_{n^{\prime}}^{\prime}\right) \in \mathbb{C}^{n^{\prime}}$ and any function $f: \mathbb{C} \rightarrow \mathbb{C}$. Here we allow the empty product and set it to be 1.

We introduce a set of contours in the complex plane as follows: Let $\Gamma_{\LL, \text {in}}, \Gamma_{\LL}$, and $\Gamma_{\LL, \text {out }}$ be three distinct contours situated within the left half-plane, each originating from $e^{-2 \pi i / 3} \infty$ and terminating at $e^{2 \pi i / 3} \infty$. Notably, $\Gamma_{\mathrm{L}, \mathrm{in}}$ corresponds to the leftmost contour, while $\Gamma_{\mathrm{L}, \text {out }}$ represents the rightmost one. The labels "in" and "out" denote their relative positions concerning $-\infty$.

\begin{figure}
    \centering
    \begin{tikzpicture}[scale=0.7]
    \draw [->,thick](-4,0) -- (4,0) node[right] {$\mathbb{R}$};
    \draw [->,thick](0,-4) -- (0,4) node[above] {$\mathrm{i}\mathbb{R}$};
    \filldraw [black] (0,0) circle (1.5pt) node[align=left, below] {$0$};
    \draw [->,orange, thick](-5,-4) -- (-3,0)--(-5,4) node[above] {$\Gamma_{\LL, \text {in}}$};
    \draw [->,red, thick](-4,-4) -- (-2,0)--(-4,4) node[above] {$\Gamma_{\LL}$};
    \draw [->,orange, thick](-3,-4) -- (-1,0)--(-3,4) node[above] {$\Gamma_{\LL, \text {out}}$};
    \draw [->,blue,thick](5,-4) -- (3,0)--(5,4) node[above] {$\Gamma_{\RR, \text {in}}$};
    \draw [->,green,thick](4,-4) -- (2,0)--(4,4) node[above] {$\Gamma_{\RR}$};
    \draw [->,blue,thick](3,-4) -- (1,0)--(3,4) node[above] {$\Gamma_{\RR, \text {out}}$};
    \end{tikzpicture}  
    \caption{Integration of contours in Proposition \ref{prop1}.}
    \label{contours}
\end{figure}
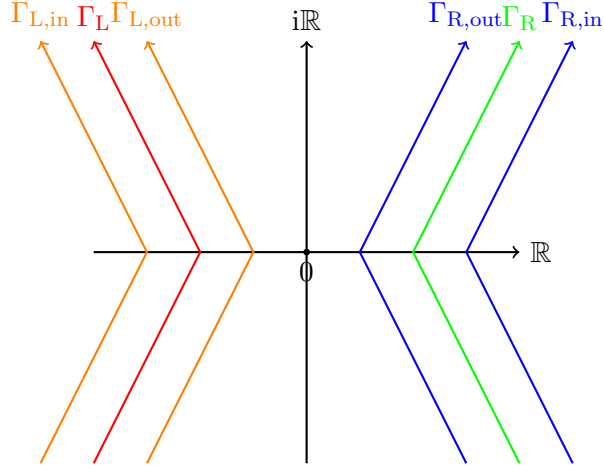

Similarly, let $\Gamma_{\mathrm{R}, \text {in }}, \Gamma_{\mathrm{R}}$, and $\Gamma_{\mathrm{R}, \text {out }}$ denote three distinct contours residing within the right half-plane. These contours originate from $e^{-\pi i / 3} \infty$ and terminate at $e^{\pi i / 3} \infty$. In this context, "in" and "out" signify their positions relative to $+\infty.$ Therefore, $\Gamma_{\mathrm{R}, \text {in }}$ represents the rightmost contour, while $\Gamma_{\mathrm{R}, \text {out }}$ corresponds to the leftmost contour. For visual clarity, refer to Figure \ref{contours} for an illustration of these contour configurations.

In the subsequent section, we derive an explicit expression for $\prob\left(\mathrm{H}(\hat{x},\hat\tau)\ge \hat\rh,\mathrm{H}(x',\tau')\ge \rh'\right)$.

\begin{prop}\label{prop1}
Suppose $\tau,\tau'>0$ and $x,x'\in \realR$ are all fixed. Then
\begin{equation}
    \begin{split}
       &\prob\left(\mathrm{H}(x',\tau')\ge \rh',\mathrm{H}(x+x',\tau+\tau')\ge \rh+\rh'\right)\\
       &=\oint_{|\rz|>1}\frac{\mathrm{d}\rz}{2\pi\mathrm{i}\rz(1-\rz)}\sum_{{n_1,n_2\ge 1}}\frac{(-1)^{n_1+n_2}}{(n_1!)^2(n_2!)^2}\rt_{n_1,n_2}(\rz;x_,x',\tau,\tau',\rh,\rh'),
    \end{split}\label{eq2}
\end{equation}
with
\begin{equation}
    \begin{split}
        &\rt_{n_1,n_2}(\rz;x_,x',\tau,\tau',\rh,\rh')\\
        &:=\prod_{i_2=1}^{n_2}\left[\frac{1}{1-\rz}\int_{\Gamma_{\LL,\inn}}\ddbar{\xi_{i_2}^{(2)}}{}-\frac{\rz}{1-\rz}\int_{\Gamma_{\LL,\out}}\ddbar{\xi_{i_2}^{(2)}}{}\right]\prod_{i_1=1}^{n_1}\int_{\Gamma_{\LL}}\ddbar{\xi_{i_1}^{(1)}}{}\\
        &\quad\prod_{i_2=1}^{n_2}\left[\frac{1}{1-\rz}\int_{\Gamma_{\RR,\inn}}\ddbar{\eta_{i_2}^{(2)}}{}-\frac{\rz}{1-\rz}\int_{\Gamma_{\RR,\out}}\ddbar{\eta_{i_2}^{(2)}}{}\right]\prod_{i_1=1}^{n_1}\int_{\Gamma_{\RR}}\ddbar{\eta_{i_1}^{(1)}}{}\\
        &\quad\left(1-\rz\right)^{n_1}\left(1-\rz^{-1}\right)^{n_2}
		F\left(\boldsymbol{\xi}^{(1)},\boldsymbol{\xi}^{(2)},\boldsymbol{\eta}^{(1)},\boldsymbol{\eta}^{(2)}\right) \label{tstar}.
    \end{split}
\end{equation}
Here
\begin{equation}
    F\left(\boldsymbol{\xi}^{(1)},\boldsymbol{\xi}^{(2)},\boldsymbol{\eta}^{(1)},\boldsymbol{\eta}^{(2)}\right):=\prod_{\ell=1}^2\mathrm{f}_\ell(\boldsymbol{\xi}^{(\ell)})\mathrm{f}_\ell(\boldsymbol{\eta}^{(\ell)})\left(C(\boldsymbol{\xi}^{(\ell)};\boldsymbol{\eta}^{(\ell)})\right)^2\cdot \frac{\Delta(\boldsymbol{\xi}^{(1)};\boldsymbol{\eta}^{(2)})}{\Delta(\boldsymbol{\xi}^{(1)};\boldsymbol{\xi}^{(2)})}\cdot \frac{\Delta(\boldsymbol{\eta}^{(1)};\boldsymbol{\xi}^{(2)})}{\Delta(\boldsymbol{\eta}^{(1)};\boldsymbol{\eta}^{(2)})},
\end{equation}
where the vectors $\boldsymbol{\xi}^{(\ell)}=\left(\xi_1^{(\ell)},\cdots,\xi_{i_\ell}^{(\ell)}\right)$ and $\boldsymbol{\eta}^{(\ell)}=\left(\eta_1^{(\ell)},\cdots,\eta_{i_\ell}^{(\ell)}\right)$ for $\ell\in\{1,2\}$ and the functions $\mathrm{f}_{\ell}$ are defined by
\begin{equation}
   \begin{split}
       &\mathrm{f}_{1}(\zeta)=
	      \begin{cases}
	     	\exp{\left(-\frac{\tau'}{3}\zeta^3+x'\zeta^2+\rh'\zeta\right)}, &\mathrm{Re}\zeta<0,\\
		    \exp{\left(\frac{\tau'}{3}\zeta^3-x'\zeta^2-\rh'\zeta\right)}, &\mathrm{Re}\zeta>0, 
	      \end{cases}\\
       &\mathrm{f}_{2}(\zeta)=
	      \begin{cases}
	     	\exp{\left(-\frac{\tau}{3}\zeta^3+x\zeta^2+\rh\zeta\right)}, &\mathrm{Re}\zeta<0,\\
		    \exp{\left(\frac{\tau}{3}\zeta^3-x\zeta^2-\rh\zeta\right)}, &\mathrm{Re}\zeta>0. \label{airy}
	      \end{cases}
   \end{split}
\end{equation}
\end{prop}

\begin{proof}
We start with
\begin{equation}
	\begin{split}
		\prob\left(\mathrm{H}(x',\tau')\ge \rh',\mathrm{H}(x+x',\tau+\tau')\ge \rh+\rh'\right)=P_1-P_2
	\end{split}
\end{equation}
where 
\begin{equation}
    \begin{split}
        &P_1:=\prob\left(\mathrm{H}(x',\tau')\ge  \rh'\right)\\
        &P_2:=\prob\left(\mathrm{H}(x',\tau')\ge \rh',\mathrm{H}(x+x',\tau+\tau')\le \rh+\rh'\right).
    \end{split}
\end{equation}
We have  
\begin{equation}
    \begin{split}
        P_1=1-\sum_{n_1\geq 0}\frac{(-1)^{n_1}}{(n_1!)^2}\prod_{i_1=1}^{n_1}\int_{\Gamma_{\LL}}\ddbar{\xi_{i_1}^{(1)}}{}\int_{\Gamma_{\RR}}\ddbar{\eta_{i_1}^{(1)}}{}\left(C(\boldsymbol{\xi}^{(1)};\boldsymbol{\eta}^{(1)})\right)^2\mathrm{f}_1(\boldsymbol{\xi}^{(1)})\mathrm{f}_1(\boldsymbol{\eta}^{(1)})
    \end{split} \label{onepoint}
\end{equation}
and by \cite[formula (2.4)]{Liu-Wang22}, 
\begin{equation}
    \begin{split}
        P_2&=-\oint_{|\rz|>1}\frac{\mathrm{d}\rz}{2\pi\mathrm{i}(1-\rz)\rz}\sum_{{n_1,n_2\ge 0}}\frac{(-1)^{n_1+n_2}}{(n_1!)^2(n_2!)^2}\rt_{n_1,n_2}(\rz;x,x',\tau,\tau',\rh,\rh'). \label{two point}
    \end{split}
\end{equation}

After evaluating the integral on $\rz$, we obtain
\begin{equation}
    \begin{split}
       P_2=P_1-\oint_{|\rz|>1}\frac{\mathrm{d}\rz}{2\pi\mathrm{i}(1-\rz)\rz}\sum_{{n_1,n_2\ge 1}}\frac{(-1)^{n_1+n_2}}{(n_1!)^2(n_2!)^2}\rt_{n_1,n_2}(\rz;x,x',\tau,\tau',\rh,\rh').
    \end{split}
\end{equation}

As a result, Proposition \ref{prop1} follows immediately.
\end{proof}

\subsection{Explicit formula for the conditional tail distribution.}

Define
\begin{equation}
    \begin{split}
        &\Tilde{\rt}_{n_1,n_2}(\rz;x,x',\tau,\tau',\rh,\rh')\\
        &:=\prod_{i_2=1}^{n_2}\left[\frac{1}{1-\rz}\int_{\Gamma_{\LL,\inn}}\ddbar{\xi_{i_2}^{(2)}}{}-\frac{\rz}{1-\rz}\int_{\Gamma_{\LL,\out}}\ddbar{\xi_{i_2}^{(2)}}{}\right]\prod_{i_1=1}^{n_1}\int_{\Gamma_{\LL}}\ddbar{\xi_{i_1}^{(1)}}{}\\
        &\quad\prod_{i_2=1}^{n_2}\left[\frac{1}{1-\rz}\int_{\Gamma_{\RR,\inn}}\ddbar{\eta_{i_2}^{(2)}}{}-\frac{\rz}{1-\rz}\int_{\Gamma_{\RR,\out}}\ddbar{\eta_{i_2}^{(2)}}{}\right]\prod_{i_1=1}^{n_1}\int_{\Gamma_{\RR}}\ddbar{\eta_{i_1}^{(1)}}{}\\
        &\quad\left(1-\rz\right)^{n_1}\left(1-\rz^{-1}\right)^{n_2}\cdot \Tilde{F}\left(\boldsymbol{\xi}^{(1)},\boldsymbol{\xi}^{(2)},\boldsymbol{\eta}^{(1)},\boldsymbol{\eta}^{(2)}\right), 
    \end{split}
\end{equation}
where
\begin{equation}
    \begin{split}
    \Tilde{F}\left(\boldsymbol{\xi}^{(1)},\boldsymbol{\xi}^{(2)},\boldsymbol{\eta}^{(1)},\boldsymbol{\eta}^{(2)}\right):=\left(\sum_{i_1=1}^{n_1}\left(\xi_{i_1}^{(1)}-\eta_{i_1}^{(1)}\right)-\sum_{i_2=1}^{n_2}\left(\xi_{i_2}^{(2)}-\eta_{i_2}^{(2)}\right)\right)F\left(\boldsymbol{\xi}^{(1)},\boldsymbol{\xi}^{(2)},\boldsymbol{\eta}^{(1)},\boldsymbol{\eta}^{(2)}\right).
    \end{split}
\end{equation}
Subsequently, we can derive an expression for the conditional tail probability.
\begin{prop} \label{explicit formula}
We have 
\begin{equation}
    \prob\left(\mathrm{H}(x+x',\tau+\tau')-\mathrm{H}(x',\tau')\ge \rh\mid \mathrm{H}(x',\tau')=\rh'\right)=\frac{\Tilde{Q}(x',x,\tau',\tau,\rh',\rh)}{-\mathrm{F}'(\rh';x',\tau')} \label{conditional2}
\end{equation}
where $\mathrm{F}'(\rh';x',\tau')$ is the derivative of $\mathrm{F}(\rh';x',\tau')$ in the $\rh'$ variable, and
\begin{equation}
    \Tilde{Q}(x',x,\tau',\tau,\rh',\rh)=\oint_{|\rz|>1}\frac{\mathrm{d}\rz}{2\pi\mathrm{i}\rz(1-\rz)}\sum_{{n_1,n_2\ge 1}}\frac{(-1)^{n_1+n_2}}{(n_1!)^2(n_2!)^2}\Tilde{\rt}_{n_1,n_2}(\rz;x',x,\tau',\tau,\rh',\rh).        \label{numerator}
\end{equation}
\end{prop}
\begin{proof}
It follows immediately from 
\begin{equation}
    \begin{split}
        \Tilde{Q}(x',x,\tau',\tau,\rh',\rh)=\frac{\partial}{\partial\rh'}\prob\left(\mathrm{H}(x',\tau')\ge \rh',\mathrm{H}(\hat{x},\hat\tau)\ge \hat\rh\right)\bigg|_{\hat{h}=h+h', \hat{x}=x+x', \hat{\tau}=\tau+\tau'}
    \end{split}
\end{equation}
and 
\begin{equation}
    -\mathrm{F}'(\rh';x',\tau')=\frac{\mathrm{d}}{\mathrm{d}\rh'}\prob\left(\mathrm{H}(x',\tau')\ge \rh'\right). \label{denominator}
\end{equation}
\end{proof}

\section{Proof of Theorem \ref{thm:main}.} \label{Main Term Section}

\section*{Outline}

The proof of Theorem \ref{thm:main} relies on a classical steepest descent analysis of the contour integral formula of the conditional probability. To provide a clear perspective, our approach commences with the derivation of this conditional probability, using the two-point joint distribution of the KPZ fixed point as a foundation. The formula, which is derived and presented in Proposition \ref{explicit formula}, takes the form of a contour integral of a series expansion resembling Fredholm-type determinants, making it conducive to asymptotic analysis. The dominant contribution to the integral arise from the terms $\Tilde{\rt}_{n_1,n_2}$ for $n_1=1, n_2\in \mathbb{N}$. The precise asymptotics coming form the $\Tilde{\rt}_{1,n_2}$ terms are proven in Section \ref{Main Term Section} and the appendix by expanding out the Vandermonde determinant expression mentioned in \eqref{vandermonde} and comparing the resulting expansion with known formulas for the Tracy-Widom distribution and its derivatives. We refer readers to the comprehensive breakdown provided in Section \ref{Main Term Section}.

However, in Lemma \ref{mainlm}, a significant challenge arises from the presence of singularities in the integrand during the deformation of contours. To address this issue, our approach initiates with an analysis in Section 3.1, where we estimate the Airy-type function along the steepest descent path. This analysis yields the primary component of the upper bound detailed in Lemma \ref{mainlm}. 

Prior to letting $\rh'$ approach infinity, we employ a switch between $\Gamma_{\LL} (\Gamma_\RR)$ and $\Gamma_{\LL,\inn} (\Gamma_{\RR,\inn})$. This strategic maneuver ensures that no further singularities are encountered during the process of letting $\rh'$ tend to infinity. To counterbalance this switch, Sections 4.2 and 4.3 require a comprehensive assessment of the iterative residue contributions resulting from this manoeuvre. We achieve this through combinatorial analysis, ultimately demonstrating that their impact is negligible when compared to the principal bound. Consequently, we can safely proceed with letting $\rh'$ tend to infinity. 
We need an estimate of the kernel $\Tilde{\rt}_{n_1,n_2}(\rz;x,x',\tau,\tau',\rh,\rh')$ to control the error term when we perform the asymptotic analysis. 
The following upper bound, which we present below, ensures the validity of taking the limit for $\rh'$ in the series. The detailed proof of this technical lemma will be provided in Section \ref{mainpf}.

We note that estimating the error terms in this paper presents greater complexity than in \cite{Liu-Wang22}. In \cite{Liu-Wang22}, the authors analyze a condition at a later time, where the primary contributions arise from cases where $n_1=n_2=1$, and errors associated with these deformed contours are more straightforward to bound. Conversely, in our study, which considers conditions at an earlier time, the two heights vary significantly, with one approaching infinity and the other remaining fixed. This discrepancy results in non-trivially ordered critical points for the multiple integrals, necessitating residue estimation when deforming contours to these critical points. Additionally, the significant contributions emerge from infinitely many terms when $n_1=1$ but $n_2$ can vary arbitrarily, corresponding to the GUE Tracy-Widom distribution. These factors combine to make the residue error estimation for $2^{n_1}$ possible scenarios for each $n_1$ and a given $n_2$ particularly challenging, as discussed in Section 4.3.

\begin{lm}\label{mainlm}
Assume $\tau,\tau',x,x'$ are all fixed positive constants and suppose $|\rz|=r>1$ is fixed. Then for any fixed $\varepsilon\in (0,1)$ and fixed $\rh,$ there exists a positive constant $C=C(r,\tau,\tau',x,x',\rh)$ such that
\begin{equation}
    \begin{split}
        &|\Tilde{\rt}_{n_1,n_2}(\rz;x,x',\tau,\tau',\rh,\rh')|\\
        &\leq C^{n_1+n_2}(n_1\land n_2)^2n_1^{1+n_1/2}n_2^{1+n_2/2}(n_1+n_2)^{(n_1+n_2)}e^{-\frac{4(1-\varepsilon)}{3}\left(n_1\tau'^{-2}(\rh'\tau'+x'^2)^{3/2}\right)} \label{lm1}
    \end{split}
\end{equation}
for positive $\rh'$ large enough.
\end{lm}

\begin{rmk}
The upper bound on $\Tilde{\rt}_{n_1,n_2}(\rz;x,x',\tau,\tau',\rh,\rh')$ enables us to manage $\Tilde{Q}(x',x,\tau',\tau,\rh',\rh)$ through a convergent series. Specifically, by applying Stirling's formula $n!\sim \sqrt{2\pi n}(n/e)^n$, we obtain
\begin{equation}
    \begin{split}
        &\sum_{n_1,n_2\ge1}\frac{(n_1\land n_2)^2}{(n_1!n_2!)^2}n_1^{1+n_1/2}n_2^{1+n_2/2}(n_1+n_2)^{(n_1+n_2)}\\
        &\leq \sum_{n_1,n_2}\frac{C^{n_1+n_2}(n_1\land n_2)^2n_1^{1+3n_1/2}n_2^{1+3n_2/2}}{n_1n_2n_1^{2n_1}n_2^{2n_2}}=\sum_{n_1,n_2}\frac{C^{n_1+n_2}(n_1\land n_2)^2}{n_1^{n_1/2}n_2^{n_2/2}}<\infty,
    \end{split}
\end{equation}
where the first inequality follows from  the following simple inequality
\begin{equation}
\left(\frac{n+m}{2}\right)^{(n+m)/2} \le n^{n/2} m^{m/2}
\end{equation}
since the function $x\ln x$ is a convex function on $(0,\infty)$.
\end{rmk}

Moreover, we can write
\begin{equation}
    \begin{split}
       &\Tilde{Q}(x',x,\tau',\tau,\rh',\rh)\\
    	&=\oint_{>1}\frac{\mathrm{d}\rz}{2\pi\mathrm{i}\rz(1-\rz)}\left(\sum_{n_1=1,n_2\ge 1}+\sum_{n_1\geq 2, n_2\ge 1}\right)\frac{(-1)^{n_1+n_2}}{(n_1!n_2!)^2}\Tilde{\rt}_{n_1,n_2}(\rz;x',x,\tau',\tau,\rh',\rh), \label{sum}
    \end{split}
\end{equation}
where  $\oint_{>1}$ is integrated over a circle around the origin in the counter-clockwise direction with radius strictly larger than 1.
By \eqref{lm1}, we obtain
\begin{equation}
    \begin{split}
        \left| \oint_{>1}\frac{\mathrm{d}\rz}{2\pi\mathrm{i}\rz(1-\rz)}\sum_{n_1\geq 2, n_2\ge 1}\frac{(-1)^{n_1+n_2}}{(n_1!n_2!)^2}\Tilde{\rt}_{n_1,n_2}(\rz;\tau,\tau',\rh,\rh')\right| \leq \Tilde{C} e^{-\frac{4(1-\varepsilon)}{3}\left(\tau'^{-2}(\rh'\tau'+x'^2)^{3/2}\right)}
    \end{split}
\end{equation}
for sufficiently large $h'$, where the constant $\Tilde{C}$ is given by 
\begin{equation}
    \Tilde{C}=\oint_{>1}\frac{|\mathrm{d}\rz|}{2\pi|\rz(1-\rz)|}\sum_{n_1\geq 2, n_2\ge 1}\frac{(n_1\land n_2)^2n_1^{1+n_1/2}n_2^{1+n_2/2}(n_1+n_2)^{(n_1+n_2)}}{(n_1!n_2!)^2} C^{n_1+n_2}<\infty.
\end{equation}

Now if we consider the case $n_1=1$ in $\Tilde{\rt}_{n_1,n_2}$ and evaluate  $\oint_{>1} \ddbar{\rz}{\rz(1-\rz)}\Tilde{\rt}_{1,n_2}$, we see that only one of the terms survives. In fact, it can be derived directly by the following computation

\begin{equation}
    \begin{split}
        \oint_{>1} (1-\rz)^{-n_2} (-\rz)^{k-n_2}\ddbar{\rz}{\rz} = 
        \begin{dcases}
            1 , & k=0 \\
            0 ,         & 1 \leq k \leq 2n_2.
        \end{dcases}
    \end{split} \label{lm2.3}
\end{equation}

So applying \eqref{lm2.3} and making the change of variables $\xi_1=\xi_1^{(1)}+x'/\tau',\eta_1=\eta_1^{(1)}+x'/\tau'$, and $\xi_{i}^{(2)}\mapsto\xi_{i}^{(2)}+x'/\tau',\eta_i^{(2)}\mapsto \eta_i^{(2)}+x'/\tau'$ for each $i$, we have
\begin{equation}
    \begin{split}
        &\oint_{|z|>1} \ddbar{\rz}{\rz(1-\rz)}\Tilde{\rt}_{1,n_2}(\rz;x,x',\tau,\tau',\rh,\rh')\\
        &=\prod_{i_2=1}^{n_2}\int_{\Gamma_{\LL,\out}-x'/\tau'}\ddbar{\xi_{i_2}^{(2)}}{}\int_{\Gamma_{\RR,\out}-x'/\tau'}\ddbar{\eta_{i_2}^{(2)}}{}\int_{\Gamma_{\LL}-x'/\tau'}\ddbar{\xi_1}{}\int_{\Gamma_{\RR}-x'/\tau'}\ddbar{\eta_1}{}\\
        &\cdot\left(\left(\xi_1-\eta_1\right)-\sum_{i_2=1}^{n_2}\left(\xi_{i_2}^{(2)}-\eta_{i_2}^{(2)}\right)\right) \left(C(\boldsymbol{\xi}^{(2)}+x'/\tau';\boldsymbol{\eta}^{(2)}+x'/\tau')\right)^2 \quad \frac{e^{-\frac{\tau'}{3}\xi_1^3+\left(\rh'+\frac{x'^2}{\tau'}\right)\xi_1}}{e^{-\frac{\tau'}{3}\eta_1^3+\left(\rh'+\frac{x'^2}{\tau'}\right)\eta_1}}\frac{1}{(\xi_1-\eta_1)^2}\\
        &\cdot\mathrm{f}_2(\boldsymbol{\xi}^{(2)}+x'/\tau')\mathrm{f}_2(\boldsymbol{\eta}^{(2)}+x'/\tau')\cdot \frac{\Delta(\boldsymbol{\xi}^{(2)}+x'/\tau';{\eta_1+x'/\tau'})\Delta(\boldsymbol{\eta}^{(2)}+x'/\tau';\xi_1+x'/\tau')}{\Delta(\boldsymbol{\xi}^{(2)}+x'/\tau';{\xi_1+x'/\tau'})\Delta(\boldsymbol{\eta}^{(2)}+x'/\tau';\eta_1+x'/\tau')}\\
        &=\prod_{i_2=1}^{n_2}\int_{\Gamma_{\LL,\out}}\ddbar{\xi_{i_2}^{(2)}}{}\int_{\Gamma_{\RR,\out}}\ddbar{\eta_{i_2}^{(2)}}{}\int_{\Gamma_{\LL}}\ddbar{\xi_1}{}\int_{\Gamma_{\RR}}\ddbar{\eta_1}{}\cdot\left(\left(\xi_1-\eta_1\right)-\sum_{i_2=1}^{n_2}\left(\xi_{i_2}^{(2)}-\eta_{i_2}^{(2)}\right)\right) \left(C(\boldsymbol{\xi}^{(2)};\boldsymbol{\eta}^{(2)})\right)^2 \\
        &\quad \frac{e^{-\frac{\tau'}{3}\xi_1^3+\left(\rh'+\frac{x'^2}{\tau'}\right)\xi_1}}{e^{-\frac{\tau'}{3}\eta_1^3+\left(\rh'+\frac{x'^2}{\tau'}\right)\eta_1}}\frac{1}{(\xi_1-\eta_1)^2}\cdot\mathrm{f}_2(\boldsymbol{\xi}^{(2)}+x'/\tau')\mathrm{f}_2(\boldsymbol{\eta}^{(2)}+x'/\tau')\cdot \frac{\Delta(\boldsymbol{\xi}^{(2)};{\eta_1})\Delta(\boldsymbol{\eta}^{(2)};\xi_1)}{\Delta(\boldsymbol{\xi}^{(2)};{\xi_1})\Delta(\boldsymbol{\eta}^{(2)};\eta_1)}
        \label{eq3}
     \end{split}
\end{equation}
Where we used the notation $\boldsymbol{\xi}^{(2)}+x'/\tau':= (\xi_1^{(2)}+x'/\tau',...,\xi_{n_2}^{(2)}+x'/\tau')$ and similarly for $\boldsymbol{\eta}^{(2)}$.
Also note the last step used Cauchy's Theorem to deform the contours back.

In order to simplify the asymptotic analysis for ($\ref{eq3}$) we introduce some new notations and two lemmas. We save the proofs for the appendix.

Define the operators
\begin{equation}
    \mathcal{I}_1[g]:= \int_{\Gamma_{\LL}}\ddbar{\xi_1}{}\int_{\Gamma_{\RR}}\ddbar{\eta_1}{} \quad \frac{e^{-\frac{\tau'}{3}\xi_1^3+\left(\rh'+\frac{x'^2}{\tau'}\right)\xi_1}}{e^{-\frac{\tau'}{3}\eta_1^3+\left(\rh'+\frac{x'^2}{\tau'}\right)\eta_1}}\frac{1}{(\xi_1-\eta_1)^2} g(\xi_1,\eta_1)
    \label{I_1 def}
\end{equation}
and
\begin{equation}
    \begin{split}
        &\mathcal{I}_2[g]\\
        &:= \prod_{i_2=1}^{n_2}\int_{\Gamma_{\LL,\out}}\ddbar{\xi_{i_2}^{(2)}}{}\int_{\Gamma_{\RR,\out}}\ddbar{\eta_{i_2}^{(2)}}{} \left(C(\boldsymbol{\xi}^{(2)};\boldsymbol{\eta}^{(2)})\right)^2\mathrm{f}_2(\boldsymbol{\xi}^{(2)}+x'/\tau')\mathrm{f}_2(\boldsymbol{\eta}^{(2)}+x'/\tau')g((\boldsymbol{\xi}^{(2)},\boldsymbol{\eta}^{(2)}))
    \end{split}
\end{equation}

and define 
\begin{equation}
    \rh^* := (\tau')^{-1}\left(\rh'+\frac{x'^2}{\tau'}\right).
\end{equation}

It is important to note throughout the following discussion that $\rh$, $x$, $\tau$, $x'$, and $\tau'$ are all held constant, while we allow $\rh'\to \infty$.

\begin{lm}\label{I1bound}
Let $a,b \in \mathbb{Z}$.
Then as $\rh' \to +\infty$
\begin{equation}
\begin{split}
    \mathcal{I}_1 [\xi_1^{-a} \eta_1^{-b}] &= \frac{(-1)^a}{16\pi\tau'(\rh^*)^{(a+b+3)/2}} e^{-\frac{4}{3} \tau' (\rh^*)^{3/2}}(1+O((\rh^*)^{-3/2}))\\
    &=(-1)^a(1-\mathrm{F}(\rh',x',\tau')) (\rh^*)^{-(a+b)/2}(1+O((\rh^*)^{-3/2})\\
    &= \frac{(-1)^a}{2}\mathrm{F'}(\rh',x',\tau')(\rh^*)^{-(a+b+1)/2}(1+O((\rh^*)^{-3/2})
\end{split}
\label{I_1 terms}
\end{equation}

and as a consequence
\begin{equation}\label{I_1[I_2[]] no difference terms}
    |\mathcal{I}_1[\mathcal{I}_2[\xi_1^{-a}\eta_1^{-b}(\boldsymbol{\xi}^{(2)})^{\alpha}(\boldsymbol{\eta}^{(2)})^{\beta}]]| \leq C^{n_2} (\rh^*)^{-(a+b+3)/2} e^{-\frac{4}{3} \tau' (\rh^*)^{3/2}}
\end{equation}
where $C$ is a constant depending only on $|\alpha|_{\infty}$, $ |\beta|_{\infty}$
\end{lm}

%\begin{rmk}
%    The proof can immediately be modified to prove
%\begin{equation}
%\begin{split}
   % &\int_{\Gamma_{\LL}-x'/\tau'}\ddbar{\xi_1}{}\int_{\Gamma_{\RR}-x'/\tau'}\ddbar{\eta_1}{} \quad \frac{e^{-\frac{\tau'}{3}\xi_1^3+\left(\rh'+\frac{x'^2}{\tau'}\right)\xi_1}}{e^{-\frac{\tau'}{3}\eta_1^3+\left(\rh'+\frac{x'^2}{\tau'}\right)\eta_1}}\frac{1}{(\xi_1-\eta_1)^2}(\xi_1+x'/\tau')^{-a} (\eta_1+x'/\tau')^{-b}\\
   % &=\mathcal{I}_1 [(\xi_1+x'/\tau')^{-a} (\eta_1+x'/\tau')^{-b}] \\
   % &= \frac{(-1)^a}{16\pi\tau'(\rh^*)^{(a+b+3)/2}} e^{-\frac{4}{3} \tau' (\rh^*)^{3/2}}(1+O((\rh^*)^{-1}))\\
  % &=(-1)^a(1-\mathrm{F}(h',x',\tau')) (\rh^*)^{-(a+b)/2}(1+O((\rh^*)^{-1})\\
   % &= \frac{(-1)^a}{2}\mathrm{F'}(h',x',\tau')(\rh^*)^{-(a+b+1)/2}(1+O((\rh^*)^{-1})
%\end{split}
%\label{I_1 terms modified}
%\end{equation}
%\end{rmk}

Going forward, for any multi-index $\textbf{k}=(k_1,...,k_{n_2}) \in \mathbb{N}^{n_2}$, we will take $|\textbf{k}|=\sum_{j=1}^{n_2} |k_j|$ to be its $\ell^1$ norm and $|\textbf{k}|_{\infty}=\max_{j=1,...,n_2} |k_j|$ to be its $\ell^{\infty}$ norm.
\begin{lm}
For $a,b \in \mathbb{N}$ and $\alpha,\beta,\gamma,\delta$ multi-indices of length $n_2$ then
\begin{equation}
    |\mathcal{I}_1[\mathcal{I}_2[\xi_1^{-a}\eta_1^{-b}(\boldsymbol{\xi}^{(2)})^{\alpha}(\boldsymbol{\eta}^{(2)})^{\beta} (\xi_1-\boldsymbol{\xi}^{(2)})^{-\gamma}(\eta_1-\boldsymbol{\eta}^{(2)})^{-\delta}]]|\leq C^{n_2} n_2^{n_2} (\rh^*)^{-\frac{a+b+|\gamma|+|\delta|}{2}} e^{-\frac{4}{3}\tau'(\rh^*)^{3/2}}
    \label{I_1[I_2[]] Error bound}
\end{equation}
for a constant $C$ depending only on $|\alpha|_{\infty},|\beta|_{\infty},|\gamma|_{\infty},|\delta|_{\infty}$.\\

\label{Error Lemma for n_1=1}
\end{lm}
As clarification for the notation above, if $\gamma=(1,2)$, then
$(\xi_1-\boldsymbol{\xi}^{(2)})^{-\gamma} = (\xi_1-\xi_1^{(2)})^{-1} (\xi_1-\xi_2^{(2)})^{-2}$. We postpone the proof of the above two lemmas in the Appendix.\ref{app}.

Next, we expand the Vandermonde term into a quadratic in $\xi_1^{-1},\eta_1^{-1}$ plus error terms as follows.
This allows us write the integrand into terms which can be fed into \eqref{I_1 terms} and \eqref{I_1[I_2[]] Error bound}.
\begin{equation} \label{vandermonde}
    \begin{split}
         &\frac{\Delta(\boldsymbol{\xi}^{(2)};{\eta_1})\Delta(\boldsymbol{\eta}^{(2)};{\xi_1})}{\Delta(\boldsymbol{\xi}^{(2)};{\xi_1})\Delta(\boldsymbol{\eta}^{(2)};{\eta_1})}
         =\prod_{i_2=1}^{n_2}\left(1+\frac{\eta_{i_2}^{(2)}-\xi_{i_2}^{(2)}}{\eta_1-\eta_{i_2}^{(2)}}\right)\cdot \left(1+\frac{\xi_{i_2}^{(2)}-\eta_{i_2}^{(2)}}{\xi_1-\xi_{i_2}^{(2)}}\right)\\
         &=\prod_{i_2=1}^{n_2}\left(1+(\eta_{i_2}^{(2)}-\xi_{i_2}^{(2)})\eta_1^{-1}\frac{1}{1-\eta_{i_2}^{(2)}\eta_1^{-1}}\right)\cdot \left(1+(\xi_{i_2}^{(2)}-\eta_{i_2}^{(2)})\xi_1^{-1}\frac{1}{1-\xi_{i_2}^{(2)}\xi_1^{-1}}\right)\\
         &=\prod_{i_2=1}^{n_2}\left(1+(\eta_{i_2}^{(2)}-\xi_{i_2}^{(2)})\eta_1^{-1}\left(\sum_{k=0}^{4}(\eta_{i_2}^{(2)}\eta_1^{-1})^k+\frac{(\eta_{i_2}^{(2)}\eta_1^{-1})^5}{1-\eta_{i_2}^{(2)}\eta_1^{-1}}\right)\right) \\
         &\cdot\left(1+(\xi_{i_2}^{(2)}-\eta_{i_2}^{(2)})\xi_1^{-1}\left(\sum_{k=0}^{4}(\xi_{i_2}^{(2)}\xi_1^{-1})^k+\frac{(\xi_{i_2}^{(2)}\xi_1^{-1})^5}{1-\xi_{i_2}^{(2)}\xi_1^{-1}}\right)\right)\\
         &= 1+\sum_{i_2=1}^{n_2}(\xi_1^{-1}-\eta_1^{-1})(\xi_{i_2}^{(2)}-\eta_{i_2}^{(2)}) +\sum_{i_2=1}^{n_2} (\xi_{i_2}^{(2)}\xi_1^{-2}-\eta_{i_2}^{(2)}\eta_1^{-2})(\xi_{i_2}^{(2)}-\eta_{i_2}^{(2)}) \\
         &+\frac{1}{2}\sum_{i_2,j_2=1}^{n_2}((\xi_1^{-2}-2\xi_1^{-1} \eta_1^{-1}+\eta_1^{-2}) - \delta_{i_2 j_2}(\xi_1^{-2}+\eta_1^{-2}))(\xi_{i_2}^{(2)}-\eta_{i_2}^{(2)})(\xi_{j_2}^{(2)}-\eta_{j_2}^{(2)}) \\
         &+E(\xi_1,\eta_1, \boldsymbol{\xi}^{(2)}, \boldsymbol{\eta}^{(2)}). 
    \end{split}
\end{equation}
Here $E(\xi_1,\eta_1, \boldsymbol{\xi}^{(2)}, \boldsymbol{\eta}^{(2)})$ is a sum of $O(n_2^2)$ terms of either the form $\xi_1^{-a}\eta_1^{-b}(\boldsymbol{\xi}^{(2)})^{\alpha}(\boldsymbol{\eta}^{(2)})^{\beta}$ for $a+b \geq 3$, $|\alpha|_{\infty}+ |\beta|_{\infty} \leq 10$, or $\xi_1^{-a}\eta_1^{-b}(\boldsymbol{\xi}^{(2)})^{\alpha}(\boldsymbol{\eta}^{(2)})^{\beta} (\xi_1-\boldsymbol{\xi}_{i_2}^{(2)})^{-\gamma}(\eta_1-\boldsymbol{\eta}_{i_2}^{(2)})^{-\delta}$
where $a+b\geq 5$, $|\alpha|_{\infty}+|\beta|_{\infty} \leq 12$, $ |\gamma|+|\delta| \geq 1$, and $ |\gamma|_{\infty}+|\delta|_{\infty} \leq 2$. Thus $(a+b+|\gamma|+|\delta|)/2 \geq 3$.

So applying (\ref{I_1 terms}) and (\ref{I_1[I_2[]] no difference terms}) to the terms in $E(\xi_1,\eta_1, \boldsymbol{\xi}^{(2)}, \boldsymbol{\eta}^{(2)})$ we see that $|\mathcal{I}_1[\mathcal{I}_2[((\xi_1-\eta_1)-\sum_{i=1}^{n_2}(\xi_i^{(2)}-\eta_i^{(2)}))E(\xi_1,\eta_1, \boldsymbol{\xi}^{(2)}, \boldsymbol{\eta}^{(2)})]| \leq C^{n_2} n_2^{n_2+2} (\rh^*)^{-3} e^{-\frac{4}{3}\tau'(\rh^*)^{3/2}}\lesssim C^{n_2} n_2^{n_2+2} (\rh^*)^{-3/2}\mathrm{F'}(\rh';x',\tau')$, and therefore by Stirling's formula,
\begin{equation}
    \begin{split}
        &\sum_{n_2=1}^{\infty}\frac{(-1)^{1+n_2}}{(n_2!)^2}\mathcal{I}_1\left[\mathcal{I}_2\left[((\xi_1-\eta_1)-\sum_{i=1}^{n_2}(\xi_i^{(2)}-\eta_i^{(2)}))E(\xi_1,\eta_1, \boldsymbol{\xi}^{(2)}, \boldsymbol{\eta}^{(2)})\right]\right] \\
        &= \bigO((\rh^*)^{-3/2} \mathrm{F'}(\rh';x',\tau'))
    \label{Error for main term}
    \end{split}
\end{equation}

Finally note that from changing variables and deforming the contour back

\begin{equation}
    \mathcal{I}_2[(\xi_i^{(2)}-\eta_i^{(2)})^k]:= \prod_{i_2=1}^{n_2}\int_{\Gamma_{\LL,\out}}\ddbar{\xi_{i_2}^{(2)}}{}\int_{\Gamma_{\RR,\out}}\ddbar{\eta_{i_2}^{(2)}}{} \left(C(\boldsymbol{\xi}^{(2)};\boldsymbol{\eta}^{(2)})\right)^2\mathrm{f}_2(\boldsymbol{\xi}^{(2)})\mathrm{f}_2(\boldsymbol{\eta}^{(2)})(\xi_i^{(2)}-\eta_i^{(2)})^k
\end{equation}

It then follows that,
\begin{equation}
    \begin{split}
        &\mathrm{F}(\rh;x,\tau) = 1+\sum_{n_2=1}^{\infty}\frac{(-1)^{n_2}}{(n_2!)^2} \mathcal{I}_2[1]\\
        &\mathrm{F'}(\rh;x,\tau) = \sum_{n_2=1}^{\infty}\frac{(-1)^{n_2}}{(n_2!)^2} \mathcal{I}_2\left[\sum_{i=1}^{n_2}(\xi_i^{(2)}-\eta_i^{(2)})\right]\\
        &\mathrm{F''}(\rh;x,\tau) = \sum_{n_2=1}^{\infty}\frac{(-1)^{n_2}}{(n_2!)^2} \mathcal{I}_2\left[\left(\sum_{i=1}^{n_2}(\xi_i^{(2)}-\eta_i^{(2)})\right)^2 \right]
    \end{split}
    \label{Derivatives of TW}
\end{equation}

Therefore, in consideration of (\ref{sum}), we can apply (\ref{I_1 terms}) to the leading terms, and exploit (\ref{Derivatives of TW}) to obtain the following asymptotic expansion,
\begin{equation}
    \begin{split}
        &\Tilde{Q}(x,x',\tau,\tau',\rh,\rh')\\
        &=\sum_{n_2=1}^{\infty}\frac{(-1)^{1+n_2}}{(n_2!)^2}\mathcal{I}_1\left[\mathcal{I}_2\left[\left((\xi_1-\eta_1)-\sum_{i=1}^{n_2}(\xi_i^{(2)}-\eta_i^{(2)})\right) \frac{\Delta(\boldsymbol{\xi}^{(2)};{\eta_1})\Delta(\boldsymbol{\eta}^{(2)};{\xi_1})}{\Delta(\boldsymbol{\xi}^{(2)};{\xi_1})\Delta(\boldsymbol{\eta}^{(2)};{\eta_1})} \right] \right]+\bigO\left(e^{-\frac{8}{3}(1-\epsilon)(\tau')^{-1/2}(\rh^*)^{3/2}}\right)\\
        &=-\mathrm{F}'(\rh';x',t')\left(1+\bigO((\rh^*)^{-3/2})\right)\left((1-\mathrm{F}(\rh;x,t))+2(\rh^*)^{-1/2}\mathrm{F}'(\rh;x,t)-2(\rh^*)^{-1}\mathrm{F}''(\rh;x,t)+\bigO((\rh^*)^{-3/2})\right)\\
        &\quad-\mathrm{F}'(\rh';x',t')\left(1+\bigO((\rh^*)^{-3/2})\right)\left(-\frac{1}{2}(\rh^*)^{-1/2}\mathrm{F}'(\rh;x,t)+(\rh^*)^{-1}\mathrm{F}''(\rh;x,t)+\bigO((\rh^*)^{-3/2})\right)\\
        &= -\mathrm{F}'\left(\rh';x',t')((1-\mathrm{F}(\rh;x,t))+\frac{3}{2}(\rh^*)^{-1/2}\mathrm{F}'(\rh;x,t)-(\rh^*)^{-1}\mathrm{F}''(\rh;x,t)+\bigO((\rh^*)^{-3/2})\right)\\
        &= -\mathrm{F}'(\rh';x',t')\left((1-\mathrm{F}(\rh;x,t)\right)+\frac{3}{2}\sqrt{\frac{\tau'}{\rh'}}\mathrm{F}'(\rh;x,t)-\frac{\tau'}{\rh'}\mathrm{F}''(\rh;x,t)+\bigO\left((\rh')^{-3/2})\right)
    \end{split}
\end{equation}
where the last equality uses $(\rh^*)^{-1/2}=(\tau'/\rh')^{1/2}(1+\bigO(1/\rh'))$.

As a result, we obtain the following
\begin{equation}
    \begin{split}
        &\prob\left(\mathrm{H}(x+x',\tau+\tau')-\mathrm{H}(x',\tau')\le \rh\mid \mathrm{H}(x',\tau')=\rh'\right) = 1-\frac{\Tilde{Q}(x,x',\tau,\tau',\rh,\rh')}{-\mathrm{F'}(\rh;x,\tau)}\\
        &=\mathrm{F}(\rh;x,\tau)-\frac{3}{2}\sqrt{\frac{\tau'}{\rh'}}\mathrm{F}'(\rh;x,t)+\frac{\tau'}{\rh'}\mathrm{F}''(\rh;x,t)+\bigO((\rh')^{-3/2}).
    \end{split}
\end{equation}
\section{Proof of Lemma \ref{mainlm}}\label{mainpf}

\subsection{Estimates on the exponents of $\mathrm{f}_\ell\left(\boldsymbol{\xi}^{(\ell)}\right),\mathrm{f}_\ell\left(\boldsymbol{\eta}^{(\ell)}\right)$.}
It is evident that the predominant contribution to $\Tilde{\rt}_{n_1,n_2}(\rz;x,x',\tau,\tau',\rh,\rh')$ arises from the functions $\mathrm{f}_\ell(\boldsymbol{\xi}^{(\ell)})$ and $\mathrm{f}_\ell(\boldsymbol{\eta}^{(\ell)})$ when $x$, $x'$, $\tau$, $\tau'$, $\rh$, and $|\rz|=r$ are held constant. We establish upper bounds by analysing the exponents associated with $\mathrm{f}_\ell(\boldsymbol{\xi}^{(\ell)}),\mathrm{f}_\ell(\boldsymbol{\eta}^{(\ell)})$.

\begin{lm}
Let $\mathrm{f}_1(\zeta),\mathrm{f}_2(\zeta)$ be defined as \eqref{airy}. Then, for all $\varepsilon\in (0,1),x,x',\tau,\tau',\rh$ are fixed, we have 
\begin{equation}
     \begin{split}
         &\left|\mathrm{f}_1\left(\xi_{i_1}^{(1)}\right)\right|\le \exp{\left(-\frac{2(1-\varepsilon)}{3}(\rh'\tau'+x'^2)^{3/2}\tau'^{-2}+g(u_{i_1}^{(1)};\tau')\right)}\text{ for all } \xi_{i_1}^{(1)} \in \Gamma_{\LL},\\
         &\left|\mathrm{f}_2\left(\xi_{i_2}^{(2)}\right)\right|\le C_{\varepsilon}(x,\rh,\tau)\exp{\left(-\frac{2(1-\varepsilon)}{3}\tau+g(u_{i_2}^{(2)};\tau)\right)}\text{ for all } \xi_{i_1}^{(1)} \in \Gamma_{\LL,\inn} \cup \Gamma_{\LL,\out},\\
         &\left|\mathrm{f}_1\left(\eta_{i_1}^{(1)}\right)\right|\le \exp{\left(-\frac{2(1-\varepsilon)}{3}(\rh'\tau'+x'^2)^{3/2}\tau'^{-2}+g(v_{i_1}^{(1)};\tau')\right)}\text{ for all }\eta_{i_1}^{(1)} \in \Gamma_{\RR}, \\
         &\left|\mathrm{f}_2\left(\eta_{i_2}^{(2)}\right)\right|\le C_{\varepsilon}(x,\rh,\tau)\exp{\left(-\frac{2(1-\varepsilon)}{3}\tau+g(v_{i_2}^{(2)};\tau)\right)}\text{ for all } \eta_{i_2}^{(2)} \in \Gamma_{\RR,\inn}\cup\Gamma_{\RR,\out} ,  \label{es1}
    \end{split} 
\end{equation}
where $g(u;\tau):=-\frac{8}{3}\tau u^3$ for $u\ge 0$ with a positive parameter $\tau$ and $C_{\varepsilon}(x,\rh,\tau)$ is some positive constant depending on $\varepsilon,x,\rh,\tau.$
\end{lm}
\begin{proof}
We first take a change of variables and deform the contours as following
\begin{equation}
    \begin{split}
        &\xi_{i_2}^{(2)}=-1+\varepsilon+e^{\pm 2\pi \ii/3}\cdot u_{i_2}^{(2)}, \text{ for } \xi_{i_2}^{(2)} \in \Gamma_{\LL,\out}=-1+\varepsilon+e^{\pm 2\pi \ii/3}\cdot \mathbb{R}_{\geq 0},\\
        &\xi_{i_1}^{(1)}=-\sqrt{\frac{\rh'}{\tau'}+\frac{x'^2}{\tau'^2}}+e^{\pm 2\pi \ii/3}\cdot u_{i_1}^{(1)}, \text{ for } \xi_{i_1}^{(1)} \in \Gamma_{\LL}=-\sqrt{\frac{\rh'}{\tau'}+\frac{x'^2}{\tau'^2}}+e^{\pm 2\pi \ii/3}\cdot \mathbb{R}_{\geq 0},\\
        &\xi_{i_2}^{(2)}=-1-\varepsilon+e^{\pm 2\pi \ii/3}\cdot u_{i_2}^{(2)}, \text{ for } \xi_{i_2}^{(2)} \in \Gamma_{\LL,\inn}=-1-\varepsilon+e^{\pm 2\pi \ii/3}\cdot \mathbb{R}_{\geq 0},\\
        &\eta_{i_2}^{(2)}=1-\varepsilon+e^{\pm \pi \ii/3}\cdot v_{i_2}^{(2)}, \text{ for } \eta_{i_2}^{(2)} \in \Gamma_{\RR,\out}=1-\varepsilon+e^{\pm \pi \ii/3}\cdot \mathbb{R}_{\geq 0},\\
        &\eta_{i_1}^{(1)}=\sqrt{\frac{\rh'}{\tau'}+\frac{x'^2}{\tau'^2}}+e^{\pm \pi \ii/3}\cdot v_{i_1}^{(1)},\text{ for } \eta_{i_1}^{(1)} \in \Gamma_{\RR}=\sqrt{\frac{\rh'}{\tau'}+\frac{x'^2}{\tau'^2}}+e^{\pm \pi \ii/3}\cdot \mathbb{R}_{\geq 0},\\
        & \eta_{i_2}^{(2)}=1+\varepsilon+e^{\pm \pi \ii/3}\cdot v_{i_2}^{(2)}, \text{ for } \eta_{i_2}^{(2)} \in \Gamma_{\RR,\inn}=1+\varepsilon+e^{\pm \pi \ii/3}\cdot \mathbb{R}_{\geq 0} \label{scaling1}
    \end{split}
\end{equation}
for any $\varepsilon\in\left(0,1\right).$

Define $m_1(\zeta)=-\frac{\tau'}{3}\zeta^3 +\left(\rh'+\frac{x'}{\tau'}\right)\zeta$ and write $\zeta=-\sqrt{\frac{\rh'}{\tau'}+\frac{x'^2}{\tau'^2}}+ (-1\pm \sqrt{3}\ii)t,$ where $t\geq 0$. We have 

\begin{equation}
    \begin{split}
        \mathrm{Re}(m_1(\zeta))&=-\frac{2}{3}(\rh'\tau'+x'^2)^{3/2}\tau'^{-2}-2t^2\left(\tau'\rh'+x'^2\right)^{1/2}-\frac{8}{3}\tau' t^3\\
        &\leq -\frac{2}{3}(\rh'\tau'+x'^2)^{3/2}\tau'^{-2}-\frac{8}{3}\tau' t^3.
    \end{split} \label{argument1}
\end{equation}
for all $\rh', t\geq 0.$

By a similar manner, define $m_2(\zeta)=-\frac{\tau}{3}\zeta^3 +\left(\rh+\frac{x^2}{\tau}\right)\zeta$ and write $\zeta=-1 + (-1\pm \sqrt{3}\ii)s,$ where $s\geq 0$. We have 
\begin{equation}
    \begin{split}
        \mathrm{Re}(m_2(\zeta+\varepsilon))
        &=-(1-\varepsilon)\left(\rh+\frac{x^2}{\tau}\right)+\frac{(1-\varepsilon)^3}{3}\tau\\
        &+\left(\tau(1-\varepsilon)^2-\left(\rh+\frac{x^2}{\tau}\right)\right)s-2\tau(1-\varepsilon)s^2-\frac{8}{3}\tau s^3\\
        &\leq \tau\left(\frac{(1-\varepsilon)^3}{3}-(1-\varepsilon)\Tilde{\rh}+p(s)\right)-\frac{8}{3}\tau s^3,
    \end{split}
\end{equation}
where $\Tilde{\rh}:=\tau^{-2}\left(\rh\tau+x^2\right)$ and $p(s):=\left(1-\Tilde{\rh}\right)s-2(1-\varepsilon)s^2.$ Clearly, if $\Tilde{\rh}>1$, $p(s)<0$ for all $s\geq 0.$ If $\Tilde{\rh}\leq 1,$ then $p(s)\leq \frac{\left(1-\Tilde{h}\right)^2}{8(1-\varepsilon)}$ by a direct calculation and hence 
\begin{equation}
    \begin{split}
         \mathrm{Re}(m_2(\zeta+\varepsilon))
         &\leq \tau\left(\frac{(1-\varepsilon)^3}{3}-(1-\varepsilon)\Tilde{\rh}+\frac{\left(1-\Tilde{h}\right)^2}{8(1-\varepsilon)}\right)-\frac{8}{3}\tau s^3\\
         &\leq -\frac{2(1-\varepsilon)}{3}\tau-\frac{8}{3}\tau s^3+(1-\varepsilon)\left(1-\Tilde{h}\right)+\frac{\left(1-\Tilde{h}\right)^2}{(1-\varepsilon)}.
    \end{split}
\end{equation}
for all $s\geq 0.$

By conducting a parallel calculation, we derive the following
\begin{equation}
    \begin{split}
        \mathrm{Re}(m_2(\zeta-\varepsilon))&=-(1+\varepsilon)\left(\rh+\frac{x^2}{\tau}\right)+\frac{(1+\varepsilon)^3}{3}\tau\\
        &+\left(\tau(1+\varepsilon)^2-\left(\rh+\frac{x^2}{\tau}\right)\right)s-2\tau(1+\varepsilon)s^2-\frac{8}{3}\tau s^3\\
        &\leq \tau\left(\frac{(1+\varepsilon)^3}{3}-(1+\varepsilon)\Tilde{\rh}+p_{\varepsilon}(s)\right)-\frac{8}{3}\tau s^3,
    \end{split}
\end{equation}
where $p_{\varepsilon}(s):=\left((1+\varepsilon)^2-\Tilde{\rh}\right)s-2(1+\varepsilon)s^2.$
Note that if $\Tilde{\rh}>(1+\varepsilon)^2, p_{\varepsilon}(s)<0$ for all $s \geq 0.$ If $\Tilde{\rh} \leq (1+\varepsilon)^2,$ then $p_{\varepsilon}(s)\leq\frac{\left((1+\varepsilon)^2-\Tilde{h}\right)^2}{8(1+\varepsilon)}$ and thus
\begin{equation}
   \begin{split}
        \mathrm{Re}(m_2(\zeta-\varepsilon))
        &\leq \tau\left(\frac{(1+\varepsilon)^3}{3}-(1+\varepsilon)\Tilde{\rh}+\frac{\left((1+\varepsilon)^2-\Tilde{h}\right)^2}{8(1+\varepsilon)}\right)-\frac{8}{3}\tau s^3\\
        &\leq -\frac{2(1-\varepsilon)}{3}\tau-\frac{8}{3}\tau s^3+(1+\varepsilon)\left(1-\Tilde{h}\right)+\frac{\left((1+\varepsilon)^2-\Tilde{h}\right)^2}{(1+\varepsilon)}
   \end{split}
\end{equation}
for all $s\geq 0$

Thus, for $\varepsilon\in (0,1)$, we have 
\begin{equation}
    \begin{split}
         \left|\mathrm{f}_2\left(\xi_{i_2}^{(2)}\right)\right|\le C_{\varepsilon}(x,\rh,\tau)\exp{\left(-\frac{2(1-\varepsilon)}{3}\tau+g(u_{i_2}^{(2)};\tau)\right)}\text{ for all } \xi_{i_1}^{(1)} \in \Gamma_{\LL,\inn} \cup \Gamma_{\LL,\out}
    \end{split}
\end{equation}

It can be verified directly that:
\begin{equation}
    \begin{split}
        \mathrm{Re}\left(-m_1\left(\eta_{i_1}^{(1)}\right)\right)&=-\frac{2}{3}\left(\rh'\tau'+x'^2\right)^{3/2}\tau'^{-2}-2\left(\rh'\tau'+x'^2\right)^{1/2}u^2-\frac{8}{3}\tau' u^3,\\
        \mathrm{Re}\left(-m_2\left(\eta_{i_2}^{(2)}\right)\right)&=-(1\pm \varepsilon)\left(\rh+\frac{x^2}{\tau}\right)+\frac{(1\pm \varepsilon)^3}{3}\tau\\
        &+\left(\tau(1\pm\varepsilon)^2+\left(\rh+\frac{x^2}{\tau}\right)\right)v-2\tau(1\pm\varepsilon)v^2-\frac{8}{3}\tau s^3
    \end{split}
\end{equation}
where $\eta_{i_1}^{(1)}=\sqrt{\frac{\rh'}{\tau'}+\frac{x'^2}{\tau'^2}}+(1\pm \sqrt{3}\ii) u, \eta_{i_2}^{(2)}=1\pm \varepsilon +(1\pm \sqrt{3}\ii) v$ for $u,v \geq 0.$

By the same argument, we can obtain upper bounds for the real parts. In summary, for any fixed $\varepsilon\in(0,1), x,x',h, \tau>0,\tau'>0$,we have
\begin{equation}
    \begin{split}
         &\mathrm{Re}\left(m_1\left(\xi_{i_1}^{(1)}\right)\right)\leq -\frac{2(1-\varepsilon)}{3}(\rh'\tau'+x'^2)^{3/2}\tau'^{-2}+g(u_{i_1}^{(1)};\tau')\text{ for all } \xi_{i_1}^{(1)} \in \Gamma_{\LL},\\
         &\mathrm{Re}\left(m_2\left(\xi_{i_2}^{(2)}\right)\right)\leq -\frac{2(1-\varepsilon)}{3}\tau+ C_{\varepsilon}(x,\rh,\tau)+g(u_{i_2}^{(2)};\tau)\text{ for all }\xi_{i_2}^{(2)} \in \Gamma_{\LL,\inn}\cup \Gamma_{\LL,\out},\\
         &\mathrm{Re}\left(-m_1\left(\eta_{i_1}^{(1)}\right)\right)\leq -\frac{2(1-\varepsilon)}{3}(\rh'\tau'+x'^2)^{3/2}\tau'^{-2}+g(v_{i_1}^{(1)};\tau')\text{ for all }\eta_{i_1}^{(1)} \in \Gamma_{\RR},\\
         &\mathrm{Re}\left(-m_2\left(\eta_{i_2}^{(2)}\right)\right)\leq -\frac{2(1-\varepsilon)}{3}\tau+ C_{\varepsilon}(x,\rh,\tau)+g(v_{i_2}^{(2)};\tau)\text{ for all }\eta_{i_2}^{(2)} \in \Gamma_{\RR,\inn}\cup \Gamma_{\RR,\out}.
    \end{split} 
\end{equation}
 Hence, by the fact that $|\mathrm{f}_\ell(\zeta+x/\tau)|=|\exp{(m_\ell(\zeta))}|\le \exp{\left(\mathrm{Re}(m_\ell(\zeta))\right)},$ we establish the desired estimates.
\end{proof}

\subsection{Estimates on determinants.}
We require an estimate for the Cauchy-type determinants in \eqref{tstar} to advance our proof. These determinants can be effectively bounded using Hadamard's inequality.

Let
\begin{equation}
    \begin{split}
        \prod_{\ell=1}^2\left(C(\boldsymbol{\xi}^{(\ell)};\boldsymbol{\eta}^{(\ell)})\right)^2 \cdot \frac{\Delta(\boldsymbol{\xi}^{(1)};\boldsymbol{\eta}^{(2)})\Delta(\boldsymbol{\eta}^{(1)};\boldsymbol{\xi}^{(2)})}{\Delta(\boldsymbol{\xi}^{(1)};\boldsymbol{\xi}^{(2)})\Delta(\boldsymbol{\eta}^{(1)};\boldsymbol{\eta}^{(2)})}
        =:B_1\cdot B_2\cdot B_3, \label{cauchy}
    \end{split}
\end{equation}
where
\begin{equation}
    \begin{split}
        B_1:=(-1)^{n_1(n_1-1)/2}C(\boldsymbol{\xi}^{(1)};\boldsymbol{\eta}^{(1)}),B_2:=(-1)^{n_2(n_2-1)/2}C(\boldsymbol{\xi}^{(2)};\boldsymbol{\eta}^{(2)}),
    \end{split}
\end{equation}
and 
\begin{equation}
    \begin{split}
        B_3:=&\det\left[
              \begin{array}{ c | c }
                    \begin{array}{ccc} & \vdots &  \\
                    & \cdots \quad \frac{1}{\xi_{i_1}^{(1)}-\eta_{j_1}^{(1)}} \quad \cdots \\ & \vdots &\end{array} &
                    \begin{array}{ccc} & \vdots &  \\ & \cdots \quad \frac{1}{\xi_{i_1}^{(1)}-\xi_{j_2}^{(2)}} \quad \cdots \\ & 
                    \vdots &\end{array} \\
                    \hline
                    \begin{array}{ccc} & \vdots & \\ & \cdots \quad \frac{1}{\eta_{i_2}^{(2)}-\eta_{j_1}^{(1)}} \quad \cdots \\ & \vdots &\end{array}& 
                    \begin{array}{ccc} & \vdots & \\ & \cdots \quad \frac{1}{\eta_{i_2}^{(2)}-\xi_{j_2}^{(2)}} \quad \cdots \\ & \vdots &\end{array}
              \end{array}
                  \right]_{\substack{1\leq i_1,j_1\leq n_1\\1\leq i_2,j_2\leq n_2}}\\
             =&(-1)^{n_1(n_1-1)/ 2+n_2(n_2+1)/2} C(\boldsymbol{\xi}^{(1)};\boldsymbol{\eta}^{(1)})C(\boldsymbol{\xi}^{(2)};\boldsymbol{\eta}^{(2)}) \cdot \frac{\Delta\left(\boldsymbol{\xi}^{(1)} ; \boldsymbol{\eta}^{(2)}\right) \Delta\left(\boldsymbol{\eta}^{(1)} ; \boldsymbol{\xi}^{(2)}\right)}{\Delta\left(\boldsymbol{\xi}^{(1)} ; \boldsymbol{\xi}^{(2)}\right) \Delta\left(\boldsymbol{\eta}^{(1)} ; \boldsymbol{\eta}^{(2)}\right)}. \label{det}
    \end{split}
\end{equation}

Then, by Lemma 3.3 in \cite{Liu22}, we obtain
\begin{equation}
    \begin{split}
        &\left|\prod_{\ell=1}^2\left(C(\boldsymbol{\xi}^{(\ell)};\boldsymbol{\eta}^{(\ell)})\right)^2 \cdot \frac{\Delta(\boldsymbol{\xi}^{(1)};\boldsymbol{\eta}^{(2)})\Delta(\boldsymbol{\eta}^{(1)};\boldsymbol{\xi}^{(2)})}{\Delta(\boldsymbol{\xi}^{(1)};\boldsymbol{\xi}^{(2)})\Delta(\boldsymbol{\eta}^{(1)};\boldsymbol{\eta}^{(2)})}\left(\sum_{i_1=1}^{n_1}\left(\xi_{i_1}^{(1)}-\eta_{i_1}^{(1)}\right)-\sum_{i_2=1}^{n_2}\left(\xi_{i_2}^{(2)}-\eta_{i_2}^{(2)}\right)\right) \right|\\
        &\leq n_1^{1+n_1/2}n_2^{1+n_2/2}(n_1+n_2)^{(n_1+n_2)/2}\prod_{i_1=1}^{n_1}\frac{1}{\mathrm{dist}\left(\xi_{i_1}^{(1)}\right)}\frac{1}{\mathrm{dist}\left(\eta_{i_1}^{(1)}\right)}\prod_{i_2=1}^{n_2}\frac{1}{\mathrm{dist}\left(\xi_{i_2}^{(2)}\right)}\frac{1}{\mathrm{dist}\left(\eta_{i_2}^{(2)}\right)},
    \end{split}
\end{equation}
where $\mathrm{dist}(\zeta)$ represents the shortest distance from the point $\zeta \in \Gamma_*$ to the contours $\Gamma_{\LL,\inn}\cup\Gamma_{\LL}\cup \Gamma_{\LL,\out}\cup\Gamma_{\RR,\inn}\cup\Gamma_{\RR}\cup\Gamma_{\RR,\out}\setminus \Gamma_*$, where $*\in \{\{\LL,\out\},\{\LL\}\,\{\LL,\inn\},\{\RR,\out\},\{\RR\},\{
\RR,\inn\}\}$.

For fixed $\varepsilon\in (0,1), x,x',\tau,\tau', \rh$, and $c\in \left(0,\frac{\varepsilon}{10}\tau^{-1}\sqrt{\rh\tau+x^2}\right).$  When $\rh'$ become large enough, $\tau'^{-1}\sqrt{\rh'\tau'+x'^2} > \tau^{-1}\sqrt{\rh\tau+x^2}+c$, so we have
\begin{align}
    \begin{split}
    \mathrm{dist}\left(\xi_{i_2}^{(2)}\right)&=\mathrm{dist}\left(\eta_{i_2}^{(2)}\right)=\frac{\sqrt{3}}{2}\cdot 2c=\sqrt{3}c,\\
     \mathrm{dist}\left(\xi_{i_1}^{(1)}\right)&=\mathrm{dist}\left(\eta_{i_1}^{(1)}\right)=\frac{\sqrt{3}}{2}\left(\tau'^{-1}\sqrt{\rh'\tau'+x'^2}- \tau^{-1}\sqrt{\rh\tau+x^2}+c\right)\geq \sqrt{3}c.
    \end{split}
\end{align} 
Hence we get
\begin{equation}
    \begin{split}
    &\left|\prod_{\ell=1}^2\left(C(\boldsymbol{\xi}^{(\ell)};\boldsymbol{\eta}^{(\ell)})\right)^2 \cdot \frac{\Delta(\boldsymbol{\xi}^{(1)};\boldsymbol{\eta}^{(2)})\Delta(\boldsymbol{\eta}^{(1)};\boldsymbol{\xi}^{(2)})}{\Delta(\boldsymbol{\xi}^{(1)};\boldsymbol{\xi}^{(2)})\Delta(\boldsymbol{\eta}^{(1)};\boldsymbol{\eta}^{(2)})}\left(\sum_{i_1=1}^{n_1}\left(\xi_{i_1}^{(1)}-\eta_{i_1}^{(1)}\right)-\sum_{i_2=1}^{n_2}\left(\xi_{i_2}^{(2)}-\eta_{i_2}^{(2)}\right)\right) \right|\\
    &\leq n_1^{1+n_1/2}n_2^{1+n_2/2}(n_1+n_2)^{(n_1+n_2)/2}\left(\sqrt{3}c\right)^{-2n_1-2n_2}. \label{es2}
    \end{split}
\end{equation}

\subsection{Estimates on iterated residues.}
However, we must acknowledge that as $\rh'$ becomes large, the contour $\Gamma_{\LL}(\Gamma_{\RR})$ will intersect with $\Gamma_{\LL,\inn}(\Gamma_{\RR,\inn})$. This intersection occurs due to the change of variables we introduce in \eqref{scaling1}. then we need to deform the contours $\Gamma_{\LL}, \Gamma_{\RR}$ in  \eqref{scaling1} such that they lie in between $\Gamma_{\LL,\inn},\Gamma_{\LL,\out}$ and $\Gamma_{\RR,\inn},\Gamma_{\RR,\out}$ respectively before we let $\rh'\to \infty$. 

In this section, our approach unfolds in two steps: First, we evaluate the residues that emerge as a consequence of contour deformation. Subsequently, we demonstrate that the magnitudes of different types of residues are uniformly upper-bounded. By quantifying the total number of residues and employing combinatorial analysis, we establish the desired upper bound on $\rt_{n_1,n_2}(\rz;x,x',\tau,\tau',\rh,\rh').$

Note that
\begin{equation}
    \begin{split}
        \frac{1}{1-\rz}\int_{\Gamma_{*,\inn}}-\frac{\rz}{1-\rz}\int_{\Gamma_{*,\out}}
        =\int_{\Gamma_{*,\out}}+\frac{1}{1-\rz}\left(\int_{\Gamma_{*,\inn}}-\int_{\Gamma_{*,\out}}\right)
    \end{split}
\end{equation}
where we omit the integrand for simplification and $*\in \{\LL,\RR\}.$

Therefore, we can rephrase $\rt_{n_1,n_2}(\rz;x,x',\tau,\tau',\rh,\rh')$ in the following way
\begin{equation}
    \begin{split}
        &\Tilde{\rt}_{n_1,n_2}(\rz;x,x',\tau,\tau',\rh,\rh')=\\
        &\quad\prod_{i_2=1}^{n_2}\left[\int_{\Gamma_{\LL,\out}}\ddbar{\xi_{i_2}^{(2)}}{}+\frac{1}{1-\rz}\left(\int_{\Gamma_{\LL,\inn}}-\int_{\Gamma_{\LL,\out}}\right)\ddbar{\xi_{i_2}^{(2)}}{}\right] \cdot\prod_{i_1=1}^{n_1}\int_{\Gamma_{\LL}}\ddbar{\xi_{i_1}^{(1)}}{}\\
        &\quad\prod_{i_2=1}^{n_2}\left[\int_{\Gamma_{\RR,\out}}\ddbar{\eta_{i_2}^{(2)}}{}+\frac{1}{1-\rz}\left(\int_{\Gamma_{\RR,\inn}}-\int_{\Gamma_{\RR,\out}}\right)\ddbar{\eta_{i_2}^{(2)}}{}\right]\cdot\prod_{i_1=1}^{n_1}\int_{\Gamma_{\RR}}\ddbar{\eta_{i_1}^{(1)}}{} \\
        &\quad\left(1-\rz\right)^{n_1}\left(1-\rz^{-1}\right)^{n_2} \Tilde{F}\left(\boldsymbol{\xi}^{(1)},\boldsymbol{\xi}^{(2)},\boldsymbol{\eta}^{(1)},\boldsymbol{\eta}^{(2)}\right) \label{newt}.
    \end{split}
\end{equation}

We will consider these residues case-by-case according to their sources.

\textbf{\textit{Case (1): No residues.}}

This case happens when we only multiply the terms $\frac{1}{2\pi \ii}\int_{\Gamma_{\LL,\out}}\mathrm{d}\xi_{i_2}^{(2)}$ and $\frac{1}{2\pi \ii}\int_{\Gamma_{\RR,\out}}\mathrm{d}\eta_{i_2}^{(2)}$.

Let
\begin{equation}
    \begin{split}
        &I_1:=\prod_{i_1=1}^{n_1}\prod_{i_2=1}^{n_2}\int_{\Gamma_{\LL,\out}}\ddbar{\xi_{i_2}^{(2)}}{}\int_{\Gamma_{\RR,\out}}\ddbar{\eta_{i_2}^{(2)}}{}\int_{\Gamma_{\LL}}\ddbar{\xi_{i_1}^{(1)}}{}\int_{\Gamma_{\RR}}\ddbar{\eta_{i_1}^{(1)}}{} \Tilde{F}\left(\boldsymbol{\xi}^{(1)},\boldsymbol{\xi}^{(2)},\boldsymbol{\eta}^{(1)},\boldsymbol{\eta}^{(2)}\right),
    \end{split}
\end{equation}

By \eqref{es1}, \eqref{es2} and omitting the constant $C_{\varepsilon}(x,\tau,\rh)$, we obtain
\begin{equation}
    \begin{split}
        |I_1|&\leq n_1^{1+n_1/2}n_2^{1+n_2/2}(n_1+n_2)^{(n_1+n_2)/2}\left(\sqrt{3}c\right)^{-2n_1-2n_2}\\ 
        &\quad \cdot e^{-\frac{4}{3}(1-\varepsilon)\left(n_1(\rh'\tau'+x'^2)^{3/2}\tau'^{-2}+n_2\tau\right)}\\
        &\quad \cdot\iiiint\limits_{[0,+\infty)^4}e^{\sum_{k=1}^2\sum_{j_k=1}^{n_k}g(u_{j_k}^{(k)};\tau_k)+g(v_{j_k}^{(k)};\tau_k)}\frac{\mathrm{d}u_{i_1}^{(1)}}{2\pi}\frac{\mathrm{d}u_{i_2}^{(2)}}{2\pi}\frac{\mathrm{d}v_{i_1}^{(1)}}{2\pi}\frac{\mathrm{d}v_{i_2}^{(2)}}{2\pi}\\
        &\le C^{n_1+n_2}n_1^{1+n_1/2}n_2^{1+n_2/2}(n_1+n_2)^{(n_1+n_2)/2}\cdot e^{-\frac{4}{3}(1-\varepsilon)\left(n_1(\rh'\tau'+x'^2)^{3/2}\tau'^{-2}\right)} \label{I1}
    \end{split}
\end{equation}
where we use the fact that $\int_{0}^{\infty}e^{-\frac{8}{3}u^3}du<\infty$ in the last inequality. 

\textbf{\textit{Case (2): Residues from $\Delta(\boldsymbol{\xi}^{(1)};\boldsymbol{\xi}^{(2)})$.}}

Define
\begin{equation}
    \begin{split}
        I_2^{(k)}&:=\prod_{i_1=1}^{n_1}\prod_{i_2=1}^{n_2}\int_{\Gamma_{\LL}}\ddbar{\xi_{i_1}^{(1)}}{}\int_{\Gamma_{\RR}}\ddbar{\eta_{i_1}^{(1)}}{}\int_{\Gamma_{\RR,\out}}\ddbar{\eta_{i_2}^{(2)}}{}\\
        &\quad\left[\left(\int_{\Gamma_{\LL,\inn}}-\int_{\Gamma_{\LL,\out}}\right)\ddbar{\xi_1^{(2)}}{}\cdots\left(\int_{\Gamma_{\LL,\inn}}-\int_{\Gamma_{\LL,\out}}\right)\ddbar{\xi_k^{(2)}}{}\right] \Tilde{F}\left(\boldsymbol{\xi}^{(1)},\boldsymbol{\xi}^{(2)},\boldsymbol{\eta}^{(1)},\boldsymbol{\eta}^{(2)}\right)
    \end{split}
\end{equation}
for $k=1,\cdots,n_1\land n_2.$

Note that $\Delta(\boldsymbol{\xi}^{(1)};\boldsymbol{\xi}^{(2)})=\prod_{i_1=1}^{n_1}\prod_{i_2}^{n_2}\left(\xi_{i_1}^{(1)}-\xi_{i_2}^{(2)}\right)$ has simple roots in $\xi_{i_2}^{(2)}=\xi_1^{(1)},\cdots,\xi_{n_1}^{(1)}$ for $i_2\in \{1,\cdots,n_2\}$. Then, for each $k,$ we can evaluate residues over $\xi_k^{(2)}$. Namely, if we only consider just one residue at $\xi_k^{(2)}=\xi_j^{(1)}$ for some $j=1,\cdots,n_1,$
\begin{equation}
    \begin{split}
        &\mathrm{Res}\left(\Tilde{F}\left(\boldsymbol{\xi}^{(1)},\boldsymbol{\xi}^{(2)},\boldsymbol{\eta}^{(1)},\boldsymbol{\eta}^{(2)}\right);\xi_k^{(2)}=\xi_j^{(1)}\right)\\
        &=e^{-\frac{\tau+\tau'}{3}\left(\xi_j^{(1)}\right)^3+\left(\rh+\rh'+\frac{x^2}{\tau}+\frac{x'^2}{\tau'}\right)\xi_j^{(1)}}\mathrm{f}_1(\boldsymbol{\xi}^{(1)})\mathrm{f}_1(\boldsymbol{\eta}^{(1)})\mathrm{f}_2(\boldsymbol{\eta}^{(2)})\prod_{i_2\neq k}^{n_2}\mathrm{f}(\xi_{i_2}^{(2)}) \cdot B_1\cdot B_2 \cdot B_3^{(j,n_1+k)}
    \end{split}
\end{equation}
where $B^{(i,j)}$ represents the minor of a determinant $B$ after deleting $i$-th row and $j$-th column. Indeed, since these residues are simple, multiplying $B_3$ by $\left(\xi_k^{(2)}-\xi_j^{(1)}\right),$ transforms the $j$-th row into a standard basis $e_{n_1+k}\in \realR^{n_1+n_2}$. Consequently, $B_3$ simplifies to the minor 
$B_3^{(j,n_1+k)}.$

Therefore, by \eqref{es1}, \eqref{es2}, we obtain
\begin{equation}
    \begin{split} &\left|\mathrm{Res}\left(\Tilde{F}\left(\boldsymbol{\xi}^{(1)},\boldsymbol{\xi}^{(2)},\boldsymbol{\eta}^{(1)},\boldsymbol{\eta}^{(2)}\right);\xi_k^{(2)}=\xi_j^{(1)}\right)\right|\\
        &\leq (n_1+n_2-1)^{(n_1+n_2-1)/2}n_1^{1+n_1/2}n_2^{1+n_2/2}C^{n_1+n_2-1}\cdot e^{-\frac{2(1-\varepsilon)}{3}\tau+g(u;\tau)}\\
        &\quad\prod_{i_2=1}^{n_2} e^{-\frac{2(1-\varepsilon)}{3}\tau+g(v_{i_2};\tau)}\prod_{i_2\neq k}^{n_2}e^{-\frac{2(1-\varepsilon)}{3}\tau+g(u_{i_2};\tau)}\\
        &\quad e^{-\frac{4(1-\varepsilon)}{3}n_1\tau'^{-2}(\rh'\tau'+x'^2)^{3/2}+\sum_{i_1=1}^{n_1}\left(g(u_{i_1};\tau')+g(v_{i_1};\tau')\right)}\\
        &\leq (n_1+n_2-1)^{(n_1+n_2-1)/2}n_1^{1+n_1/2}n_2^{1+n_2/2}C^{n_1+n_2-1} e^{-\frac{4(1-\varepsilon)}{3}\left(n_1\tau'^{-2}(\rh'\tau'+x'^2)^{3/2}\right)}\\
        &\quad e^{\sum_{j_1=1}^{n_1}\left(g(u_{j_1};\tau')+g(v_{j_1};\tau')\right)+\sum_{j_2=1}^{n_2}\left(g(u_{j_2};\tau)+g(v_{j_2};\tau)\right)}.  \label{oneresidue}
    \end{split}
\end{equation}

Hence, we obtain an estimate on integrals with exactly one residue at $\xi_k^{(2)}=\xi_j^{(1)}$ for some $j=1,\cdots,n_1,$
\begin{equation}
    \begin{split}
         &\left|\int_{\Gamma_{\RR,\out}}\ddbar{\eta_{i_2}^{(2)}}{}\int_{\Gamma_{\LL}}\ddbar{\xi_{i_1}^{(1)}}{}\int_{\Gamma_{\RR}}\ddbar{\eta_{i_1}^{(1)}}{}\mathrm{Res}\left(\Tilde{F}\left(\boldsymbol{\xi}^{(1)},\boldsymbol{\xi}^{(2)},\boldsymbol{\eta}^{(1)},\boldsymbol{\eta}^{(2)}\right);\xi_k^{(2)}=\xi_j^{(1)}\right)\right|\\
         &\le C^{n_1+n_2}(n_1+n_2-1)^{(n_1+n_2-1)/2}n_1^{1+n_1/2}n_2^{1+n_2/2}e^{-\frac{4}{3}(1-\varepsilon)\left(n_1(\rh'\tau'+x'^2)^{3/2}\tau'^{-2}\right)}.
    \end{split}
\end{equation}

Note $k$ represents the number of residues evaluated from $\Delta(\boldsymbol{\xi}^{(1)};\boldsymbol{\xi}^{(2)})$. Thus, for each $k$, the amount of integrals with exactly $k$ residues equals
\begin{equation}
    \binom{n_2}{k}n_1(n_1-1)\cdots(n_1-k+1)=k!\binom{n_1}{k}\binom{n_2}{k}. \label{counting}
\end{equation}

On the other hand, the exponential factor $e^{-\frac{4}{3}(1-\varepsilon)\left(n_1(\rh'\tau'+x'^2)^{3/2}\tau'^{-2}\right)}$, as a result of evaluation of residue at $\xi_k^{(2)}=\xi_j^{(1)}$, is decreasing in $n_1$. To simplify our analysis, it is advisable to disregard any additional exponential factors that arise during the evaluation of iterated residues up to the $k$ level. Therefore, by \eqref{oneresidue}, for each $k$, we have 
\begin{equation}
    \begin{split}
        \left|I_2^{(k)}\right|\le C^{n_1+n_2}k!\binom{n_1}{k}\binom{n_2}{k}(n_1+n_2-k)^{(n_1+n_2-k)/2}n_1^{1+n_1/2}n_2^{1+n_2/2}e^{-\frac{4}{3}(1-\varepsilon)\left(n_1(\rh'\tau'+x'^2)^{3/2}\tau'^{-2}\right)}.
    \end{split}
\end{equation}
The factor $(n_1+n_2-k)^{(n_1+n_2-k)/2}$ arises from applying Hadamard's type estimate to the minor of $B_3$. This estimate is computed after removing the row and column corresponding to $\left(\xi_j^{(1)}-\xi_k^{(2)}\right)^{-1}$ for a total of $k$ iterations.

If we define $b: \{k \in \mathbb{Z}_{>0}: k \leq n_1 \wedge n_2\} \to \mathbb{R}$ by
\begin{equation}
    b(k)=k!\binom{n_1}{k}\binom{n_2}{k}(n_1+n_2-k)^{(n_1+n_2-k)/2}, \label{b_k}
\end{equation}
then by the trivial bounds $k!\leq k^k$, $\binom{n_1}{k} \leq 2^{n_1}$, and $\binom{n_2}{k}\leq 2^{n_2}$ (the latter two coming from the binomial theorem), we have
\begin{equation}
    \begin{split}
        b(k)&=k!\binom{n_1}{k}\binom{n_2}{k}(n_1+n_2-k)^{(n_1+n_2-k)/2} \\
        &\leq  2^{n_1+n_2} k^k (n_1+n_2-k)^{(n_1+n_2-k)/2}\\
        &=2^{n_1+n_2}\left(\frac{k}{\sqrt{n_1+n_2}}\right)^{k}\left(1-\frac{k}{n_1+n_2}\right)^{(n_1+n_2-k)/2} (n_1+n_2)^{\frac{n_1+n_2}{2}}\\
        &\leq 2^{n_1+n_2}(n_1+n_2)^{n_1+n_2}
    \end{split} 
\end{equation}

As a result, we derive the cumulative contribution of all residues from $\Delta(\boldsymbol{\xi}^{(1)};\boldsymbol{\xi}^{(2)})$.
\begin{equation}
    \begin{split}
        \sum_{k=1}^{n_1\land n_2}\left|I_2^{(k)}\right|&\le C^{n_1+n_2}n_1^{1+n_1/2}n_2^{1+n_2/2}e^{-\frac{4}{3}(1-\varepsilon)\left(n_1(\rh'\tau'+x'^2)^{3/2}\tau'^{-2}\right)}\\
        &\quad\cdot \sum_{k=1}^{n_1\land n_2}k!\binom{n_1}{k}\binom{n_2}{k}(n_1+n_2-k)^{(n_1+n_2-k)/2}\\
        & \le C^{n_1+n_2}(n_1\land n_2)n_1^{1+n_1/2}n_2^{1+n_2/2}(n_1+n_2)^{(n_1+n_2)}e^{-\frac{4(1-\varepsilon)}{3}\left(n_1(\rh'\tau'+x'^2)^{3/2}\tau'^{-2}\right)}.
    \end{split}
\end{equation}

\textbf{\textit{Case (3): Residues from $\Delta(\boldsymbol{\eta}^{(1)};\boldsymbol{\eta}^{(2)})$.}}

Due to the symmetry of integrals, we can apply the same argument for the residues originating from $\Delta(\boldsymbol{\xi}^{(1)};\boldsymbol{\xi}^{(2)})$. In other words, we define
\begin{equation}
    \begin{split}
        I_3^{(k)}&:=\prod_{i_1=1}^{n_1}\prod_{i_2=1}^{n_2}\int_{\Gamma_{\LL}}\ddbar{\xi_{i_1}^{(1)}}{}\int_{\Gamma_{\RR}}\ddbar{\eta_{i_1}^{(1)}}{}\int_{\Gamma_{\LL,\out}}\ddbar{\xi_{i_2}^{(2)}}{}\\
        &\quad\left[\left(\int_{\Gamma_{\RR,\inn}}-\int_{\Gamma_{\RR,\out}}\right)\ddbar{\eta_1^{(2)}}{}\cdots\left(\int_{\Gamma_{\RR,\inn}}-\int_{\Gamma_{\RR,\out}}\right)\ddbar{\eta_k^{(2)}}{}\right] \Tilde{F}\left(\boldsymbol{\xi}^{(1)},\boldsymbol{\xi}^{(2)},\boldsymbol{\eta}^{(1)},\boldsymbol{\eta}^{(2)}\right) \label{I3}
    \end{split}
\end{equation}
for $k=1,\cdots,n_1\land n_2.$ We have 
\begin{equation}
    \begin{split}
        \sum_{k=1}^{n_1\land n_2}\left|I_3^{(k)}\right|&\le C^{n_1+n_2}(n_1\land n_2)n_1^{1+n_1/2}n_2^{1+n_2/2}(n_1+n_2)^{(n_1+n_2)}e^{-\frac{4(1-\varepsilon)}{3}\left(n_1(\rh'\tau'+x'^2)^{3/2}\tau'^{-2}\right)}.
    \end{split}
\end{equation}

\textbf{\textit{Case (4): Residues from $\Delta(\boldsymbol{\xi}^{(1)};\boldsymbol{\xi}^{(2)})\Delta(\boldsymbol{\eta}^{(1)};\boldsymbol{\eta}^{(2)})$. }}

Here, we need to consider contributions from iterated residues at $\xi_k^{(2)}=\xi_j^{(1)}$ and $\eta_\ell^{(2)}=\eta_m^{(1)}$ for some $j,m\in \{1,\cdots,n_1\}$. We define
\begin{equation}
    \begin{split}
        &I_4^{(k,\ell)}:=\prod_{i_1=1}^{n_1}\prod_{i_2=1}^{n_2}\int_{\Gamma_{\LL}}\ddbar{\xi_{i_1}^{(1)}}{}\int_{\Gamma_{\RR}}\ddbar{\eta_{i_1}^{(1)}}{}\left[\left(\int_{\Gamma_{\LL,\inn}}-\int_{\Gamma_{\LL,\out}}\right)\ddbar{\xi_1^{(2)}}{}\cdots\left(\int_{\Gamma_{\LL,\inn}}-\int_{\Gamma_{\LL,\out}}\right)\ddbar{\xi_k^{(2)}}{}\right]\\
        &\quad\left[\left(\int_{\Gamma_{\RR,\inn}}-\int_{\Gamma_{\RR,\out}}\right)\ddbar{\eta_1^{(2)}}{}\cdots\left(\int_{\Gamma_{\RR,\inn}}-\int_{\Gamma_{\RR,\out}}\right)\ddbar{\eta_\ell^{(2)}}{}\right] \Tilde{F}\left(\boldsymbol{\xi}^{(1)},\boldsymbol{\xi}^{(2)},\boldsymbol{\eta}^{(1)},\boldsymbol{\eta}^{(2)}\right) .\\
    \end{split}
\end{equation}
for $k,\ell=1,\cdots,n_1\land n_2.$

The order of evaluation of these iterated integrals does not impact their values, thanks to both the algebraic structure of the integrand and the deformation of contours. If we begin by evaluating the residue at $\xi_k^{(2)}=\xi_j^{(1)}$, as a consequence of \eqref{oneresidue}, we obtain:
\begin{equation}
    \begin{split}
    &\mathrm{Res}\left(\mathrm{Res}\left(\Tilde{F}\left(\boldsymbol{\xi}^{(1)},\boldsymbol{\xi}^{(2)},\boldsymbol{\eta}^{(1)},\boldsymbol{\eta}^{(2)}\right);\xi_k^{(2)}=\xi_j^{(1)}\right);\eta_\ell^{(2)}=\eta_m^{(1)}\right)\\
       &=e^{-\frac{\tau+\tau'}{3}\left(\left(\xi_j^{(1)}\right)^3+\left(\eta_m^{(1)}\right)^3\right)+\left(\rh+\rh'+\frac{x'^2}{\tau'}+\frac{x^2}{\tau}\right)\left(\xi_j^{(1)}+\eta_m^{(1)}\right)}\\
       &\quad\cdot\mathrm{f}_1(\boldsymbol{\xi}^{(1)})\mathrm{f}_1(\boldsymbol{\eta}^{(1)})\prod_{i_2\neq k}^{n_2}\mathrm{f}(\xi_{i_2}^{(2)})\prod_{i_2\neq \ell}^{n_2}\mathrm{f}_2(\eta^{(2)})\cdot B_1\cdot B_2 \cdot \left(B_3^{(j,n_1+k)}\right)^{(n_1-1+\ell,m)}. \label{I2}
    \end{split}
\end{equation}
This implies
\begin{equation}
    \begin{split}
    &\left|\mathrm{Res}\left(\mathrm{Res}\left(\Tilde{F}\left(\boldsymbol{\xi}^{(1)},\boldsymbol{\xi}^{(2)},\boldsymbol{\eta}^{(1)},\boldsymbol{\eta}^{(2)}\right);\xi_k^{(2)}=\xi_j^{(1)}\right);\eta_\ell^{(2)}=\eta_m^{(1)}\right)\right|\\
        &\le (n_1+n_2-2)^{(n_1+n_2-2)/2}n_1^{1+n_1/2}n_2^{1+n_2/2}\left(\sqrt{3}c\right)^{-2(n_1+n_2-2)}\\
        &\quad\cdot e^{-\frac{4(1-\varepsilon)}{3}\left(n_1\tau'^{-2}(\rh'\tau'+x'^2)^{3/2}+n_2\tau\right)}\cdot e^{\sum_{j_1=1}^{n_1}\left(g(u_{j_1};\tau')+g(v_{j_1};\tau')\right)}.
    \end{split}
\end{equation}
This leads us to establish a uniform upper bound for integrals involving iterated residues, both at $\xi_k^{(2)}=\xi_j^{(1)}$ and $\eta_\ell^{(2)}=\eta_m^{(1)}$.
\begin{equation}
    \begin{split}
        &\left|\int_{\Gamma_{\LL}}\ddbar{\xi_{i_1}^{(1)}}{}\int_{\Gamma_{\RR}}\ddbar{\eta_{i_1}^{(1)}}{}\mathrm{Res}\left(\mathrm{Res}\left(\Tilde{F}\left(\boldsymbol{\xi}^{(1)},\boldsymbol{\xi}^{(2)},\boldsymbol{\eta}^{(1)},\boldsymbol{\eta}^{(2)}\right);\xi_k^{(2)}=\xi_j^{(1)}\right);\eta_\ell^{(2)}=\eta_m^{(1)}\right)\right|\\
        &\le C^{n_1+n_2}(n_1+n_2-2)^{(n_1+n_2-2)/2}n_1^{1+n_1/2}n_2^{1+n_2/2}e^{-\frac{4}{3}(1-\varepsilon)\left(n_1(\rh'\tau'+x'^2)^{3/2}\tau'^{-2}\right)}.
    \end{split}
\end{equation}

We need to determine the number of iterated residues originating from $\Delta(\boldsymbol{\xi}^{(1)};\boldsymbol{\xi}^{(2)})\Delta(\boldsymbol{\eta}^{(1)};\boldsymbol{\eta}^{(2)})$. 
By the same argument as \eqref{counting}, the number of iterated residues equals
\begin{equation}
    \begin{split}
        &\sum_{k,\ell=1}^{n_1\land n_2}k!\ell!\binom{n_1}{k}\binom{n_2}{k}\binom{n_1}{\ell}\binom{n_2}{\ell}.
    \end{split}
\end{equation}
Therefore, for each $k,\ell$, we have 
\begin{equation}
    \begin{split}
         \left|I_4^{(k,\ell)}\right|&\le C^{n_1+n_2}n_1^{n_1/2}n_2^{n_2/2}e^{-\frac{4}{3}(1-\varepsilon)\left(n_1(\rh'\tau'+x'^2)^{3/2}\tau'^{-2}\right)}\\
        &\quad\cdot k!\ell!\binom{n_1}{k}\binom{n_2}{k}\binom{n_1}{\ell}\binom{n_2}{\ell}(n_1+n_2-k-\ell)^{(n_1+n_2-k-\ell)/2}.
    \end{split}
\end{equation}
The factor $(n_1+n_2-k-\ell)^{(n_1+n_2-k-\ell)/2}$ arises from applying a Hadamard-type estimate to the minor of $B_3$, wherein we remove rows and columns corresponding to $\left(\xi_j^{(1)}-\xi_k^{(2)}\right)^{-1}$ for $k$ instances and $\left(\eta_m^{(1)}-\eta_\ell^{(2)}\right)^{-1}$ for $\ell$ instances.

Using the elementary inequalities $\binom{n}{k} \leq 2^n$ and $k! \ell! \leq (k+\ell)!$, we can have
\begin{equation}
    \begin{split}
        & k!\ell!\binom{n_1}{k}\binom{n_2}{k}\binom{n_1}{\ell}\binom{n_2}{\ell}(n_1+n_2-k-\ell)^{(n_1+n_2-k-\ell)/2} \\
        &\leq  2^{2(n_1+n_2)}(k+\ell)!(n_1+n_2-k-\ell)^{(n_1+n_2-k-\ell)/2}\\
        &\leq 2^{2(n_1+n_2)}(n_1+n_2)^{(n_1+n_2)/2}\left(\dfrac{k+\ell}{\sqrt{n_1+n_2}}\right)^{k+\ell}\\
        &\leq 2^{2(n_1+n_2)}(n_1+n_2)^{(n_1+n_2)} .
    \end{split}
\end{equation}

Therefore, we obtain
\begin{equation}
    \begin{split}
        \sum_{k,\ell=1}^{n_1\land n_2}\left|I_4^{(k,\ell)}\right|\leq  C^{n_1+n_2}(n_1\land n_2)^2n_1^{1+n_1/2}n_2^{1+n_2/2}(n_1+n_2)^{(n_1+n_2)}e^{-\frac{4}{3}(1-\varepsilon)\left(n_1(\rh'\tau'+x'^2)^{3/2}\tau'^{-2}\right)}.
        \label{I4}
    \end{split}
\end{equation}

By consolidating the estimates derived in the four cases from equations \eqref{I1}, \eqref{I2}, \eqref{I3}, and \eqref{I4}, we arrive at:
\begin{equation}
    \begin{split}
         &\left|\Tilde{\rt}_{n_1,n_2}(\rz;x,x',\tau,\tau',\rh,\rh')\right|\\
         &\leq |1-\rz|^{n_1}\left|1-\rz^{-1}\right|^{n_2}\left(|I_1|+\sum_{k=1}^{n_1\land n_2}|1-\rz|^{-k}\left(\left|I_2^{(k)}\right|+\left|I_3^{(k)}\right|\right)+\sum_{k,\ell=1}^{n_1\land n_2}|1-\rz|^{-(k+\ell)}\left|I_4^{(k,\ell)}\right|\right)\\ 
         & \leq C^{n_1+n_2}(n_1\land n_2)^2n_1^{1+n_1/2}n_2^{1+n_2/2}(n_1+n_2)^{(n_1+n_2)}e^{-\frac{4}{3}(1-\varepsilon)\left(n_1(\rh'\tau'+x'^2)^{3/2}\tau'^{-2}\right)} \label{ineq}
    \end{split}
\end{equation}
for a constant $C$ independent of $\rh',$ but containing $\rh,x,\tau,c$ and $|\rz|=r>1$. This concludes the proof of Lemma \ref{mainlm}.

\bibliographystyle{alpha}
\bibliography{reference}
\Addresses
\newpage
\section{Appendix}\label{app}

Recall that $\rh^*:=(\tau')^{-1}\left(\rh'+\frac{(x')^2}{\tau'}\right)$.
\subsection{Proof of Lemma \ref{I1bound}} 
\begin{proof}
The goal of \eqref{I_1 terms} is to obtain $\rh' \to \infty$ asymptotics for the integral 

$\mathcal{I}_1 [\xi_1^{-a} \eta_1^{-b}]=\int_{\Gamma_{\LL}}\ddbar{\xi_1}{}\int_{\Gamma_{\RR}}\ddbar{\eta_1}{} \quad \frac{e^{-\frac{\tau'}{3}\xi_1^3+\left(\rh'+\frac{x'^2}{\tau'}\right)\xi_1}}{e^{-\frac{\tau'}{3}\eta_1^3+\left(\rh'+\frac{x'^2}{\tau'}\right)\eta_1}}\frac{1}{(\xi_1-\eta_1)^2} \xi_1^{-a} \eta_1^{-b}$ for $a,b \in \mathbb{Z}$

 Towards this goal we start by changing variables $\xi_1 = (\rh^*)^{1/2}z$, $\dif{\xi_1} = (\rh^*)^{1/2}\dif{z}$, and $\eta_1 = (\rh^*)^{1/2}w$, $\dif{\eta_1} = (\rh^*)^{1/2}\dif{w}$, we have
\begin{equation}
\begin{split}
    \mathcal{I}_1 [\xi_1^{-a} \eta_1^{-b}] &= -\frac{1}{4\pi^2 (\rh^*)^{(a+b)/2}} \int_{\Gamma_{\LL}^{SD}} \int_{\Gamma_{\RR}^{SD}} \exp(\tau'(\rh^*)^{3/2}((w^3-z^3)/3-(w-z))) \frac{1}{(w-z)^2} z^{-a} w^{-b} \dif{w} \dif{z}
\end{split}
\end{equation}
where after deforming the contour through Cauchy's theorem we can take $\Gamma_{\LL}^{SD}$ as the steepest descent contour passing through $z=-1$ and $\Gamma_R^{SD}$ as the steepest descent contour passing through $w=1$. If $z=\alpha+i\beta$ and $w=\gamma+i\delta$, then we can describe these contours explicitly as $\Gamma_{L}^{SD}=\{(\alpha,\beta) \in \mathbb{R}^2: \alpha = -\sqrt{1+\frac{1}{3}\beta^2}\}$, $\Gamma_{R}^{SD}=\{(\gamma,\delta) \in \mathbb{R}^2:\gamma = \sqrt{1+\frac{1}{3}\delta^2}\}$. Parameterizing these contours and applying Laplace's method we can compute

\begin{equation}
\begin{split}
    \mathcal{I}_1 [\xi_1^{-a} \eta_1^{-b}] &= \frac{1}{4\pi^2 (\rh^*)^{(a+b)/2}} \int_{\mathbb{R}} \int_{\mathbb{R}} \exp(\tau'(\rh^*)^{3/2}(-4/3-x^2-y^2)) \\
    &\cdot\left(\frac{1}{2^2}(-1)^a (1)^b +\bigO(x^2+y^2)\right)\dif{x} \dif{y} \\
    &=\frac{(-1)^a}{16\pi^2} \cdot \frac{\pi}{\tau'(\rh^*)^{(a+b+3)/2}} e^{-\frac{4}{3} \tau' (\rh^*)^{3/2}}(1+O((\rh^*)^{-3/2}))
\end{split}
\end{equation}
This proves the first equality in (\ref{I_1 terms}). The next two equalities follow from asymptotics for $\mathrm{F}_{\mathrm{GUE}}$ and $\mathrm{F'}_{\mathrm{GUE}}$, (\ref{TWasymp1}) and (\ref{TWasymp2}).

Finally to prove \eqref{I_1[I_2[]] no difference terms}, note that $\mathcal{I}_1[\mathcal{I}_2[\xi_1^{-a} \eta_1^{-b}(\boldsymbol{\xi}^{(2)})^{\alpha}(\boldsymbol{\eta}^{(2)})^{\beta}]] =\mathcal{I}_1[\xi_1^{-a} \eta_1^{-b}]\mathcal{I}_2[(\boldsymbol{\xi}^{(2)})^{\alpha}(\boldsymbol{\eta}^{(2)})^{\beta}]$ so after applying \eqref{I_1 terms}, it suffices to prove that $|\mathcal{I}_2[(\boldsymbol{\xi}^{(2)})^{\alpha}(\boldsymbol{\eta}^{(2)})^{\beta}]| \leq C(|\alpha|_{\infty},|\beta|_{\infty})^{n_2}$. For this deform the contours $\Gamma_{\LL,\out}$ and $\Gamma_{\RR,\out}$ to the steepest descent contours $\Gamma_{L}^{SD}=\{(\alpha,\beta) \in \mathbb{R}^2: \alpha = -\sqrt{1+\frac{1}{3}\beta^2}\}$, and $ \Gamma_{R}^{SD}=\{(\gamma,\delta) \in \mathbb{R}^2:\gamma = \sqrt{1+\frac{1}{3}\delta^2}\}$ respectively. 

On these contours observe that $\mathrm{f}_2(\xi^{(2)}_i+x'/\tau')\mathrm{f}_2(\eta^{(2)}_i+x'/\tau') \leq C e^{-\frac{1}{4}\tau(|\xi^{(2)}_i|^3+|\eta^{(2)}_i|^3)}$ (where the constant $C$ may depend on the parameters $\tau,\tau',x,x',\rh$), and thus

\begin{equation}
    \begin{split}
        |\mathcal{I}_2[(\boldsymbol{\xi}^{(2)})^{\alpha}(\boldsymbol{\eta}^{(2)})^{\beta}]| \leq C(|\alpha|_{\infty},|\beta|&\leq \prod_{i=1}^{n_2}\frac{C}{4\pi^2 }  \int_{\Gamma_{\LL}^{SD}} \int_{\Gamma_{\RR}^{SD}} e^{-\frac{1}{4}\tau(|\xi^{(2)}_i|^3+|\eta^{(2)}_i|^3)}|\xi^{(2)}_i|^{\alpha_i} |\eta^{(2)}_i|^{\beta_i} |\dif{w_i}| |\dif{z_i}|\\&
        \leq \prod_{i=1}^{n_2}\frac{C}{4\pi^2 }  \int_{\Gamma_{\LL}^{SD}} \int_{\Gamma_{\RR}^{SD}} e^{-\frac{1}{4}\tau(|\xi^{(2)}_i|^3+|\eta^{(2)}_i|^3)}|\xi^{(2)}_i|^{|\alpha|_{\infty}} |\eta^{(2)}_i|^{|\beta|_{\infty}} |\dif{w_i}| |\dif{z_i}|\\
        &\leq C(|\alpha|_{\infty},|\beta|_{\infty})^{n_2}
    \end{split}
    \label{I_2 of monomials bound}
\end{equation}
\end{proof}

\subsection{Proof of Lemma \ref{Error Lemma for n_1=1}}

 We first need a little estimate on the distance between points on the scaled contours.

\begin{lm}

Suppose $z,z_1\in \Gamma_{\RR}^{SD}$, the steepest descent contours,

\begin{equation}
    |(\rh^*)^{1/2} z - z_1| \geq C \frac{(\rh^*)^{1/2}}{|z_1|} \label{distance between contours}
\end{equation}
for $\rh^* \geq 4$ and some constant $C>0$. 
\label{contour difference lemma}
\end{lm}

\begin{rmk}

Of course the same conclusion holds if we replace $z$ and $z_1$ with $w,w_1 \in \Gamma_{\LL}^{SD}$
\end{rmk}

\begin{proof}

Let $z=x+iy = \sqrt{1+\frac{1}{3}y^2}+iy$ and $z_1=x_1+iy_1 = \sqrt{1+\frac{1}{3}y_1^2}+iy_1$. \\
If $|y_1| \leq 1$, then $x_1 \leq \sqrt{4/3}$, so $|(\rh^*)^{1/2} x - x_1| \geq (\rh^*)^{1/2}-\sqrt{4/3}\geq (1-\frac{1}{2}\sqrt{4/3})(\rh^*)^{1/2} \geq (1-\frac{1}{2}\sqrt{4/3})\frac{(\rh^*)^{1/2}}{|z_1|}$, so from now on suppose $|y_1| \geq 1$. \\

For now, let's also assume that $|y|\geq 1$. \\
If additionally $\left|\frac{(\rh^*)^{1/2}y}{y_1}\right| \leq 1/2$, then
$|(\rh^*)^{1/2} y-y_1| =|y_1|\cdot\left|\frac{(\rh^*)^{1/2} y}{y_1}-1 \right| 
\geq \frac{1}{2}|y_1| \geq |(\rh^*)^{1/2} y| \geq (\rh^*)^{1/2} \geq \frac{(\rh^*)^{1/2}}{|z_1|}$. \\
Conversely if $\left|\frac{y_1}{(\rh^*)^{1/2} y}\right| < 1/2$, then 
$|(\rh^*)^{1/2} y-y_1| =|(\rh^*)^{1/2} y|\cdot\left|1-\frac{y_1}{(\rh^*)^{1/2} y} \right| \geq \frac{1}{2}|(\rh^*)^{1/2} y| \geq \frac{1}{2} (\rh^*)^{1/2} \geq \frac{(\rh^*)^{1/2}}{2|z_1|}$.\\

Now we can assume $\frac{1}{2} \leq \left|\frac{(\rh^*)^{1/2} y}{y_1}\right| \leq 2$. If  $y_1^2-\rh^* y^2 \geq \rh^*$, then $|(\rh^*)^{1/2} y-y_1| \geq \frac{\rh^*}{|(\rh^*)^{1/2} y|+|y_1|} \geq \frac{\rh^*}{3|y_1|} \geq \frac{(\rh^*)^{1/2}}{3|y_1|} \geq \frac{(\rh^*)^{1/2}}{3 |z_1|} \geq C \frac{(\rh^*)^{1/2}}{|z_1|} $. \\

On the other hand if $y_1^2-\rh^*y^2 < \rh^*$

\begin{equation}
    \begin{split}
        & \left|(\rh^*)^{1/2}\sqrt{1+\frac{1}{3}y^2}-\sqrt{1+\frac{1}{3}y_1^2}\right|=\frac{|\rh^*(1+\frac{1}{3}y^2)-(1+\frac{1}{3}y_1^2)|}{(\rh^*)^{1/2}\sqrt{1+\frac{1}{3}y^2}+\sqrt{1+\frac{1}{3}y_1^2}} \\
        &\geq \frac{|\rh^*+\frac{1}{3}(\rh^* y^2- y_1^2)|}{\frac{4}{3}(|(\rh^*)^{1/2} y|+| y_1|)} \geq \frac{\frac{2}{3}\rh^*-1}{4| y_1|}\\
        &\geq C \frac{\rh^*}{ |z_1|} \geq C \frac{(\rh^*)^{1/2}}{|z_1|}
    \end{split}
\end{equation}
and hence $|(\rh^*)^{1/2} z- z_1| \geq C \frac{(\rh^*)^{1/2 }}{|z_1|}$. \\

Finally if $|y|< 1$ then either $|(\rh^*)^{1/2} z - z_1| \geq \frac{1}{10}(\rh^*)^{1/2}$ and we're done, or $|(\rh^*)^{1/2} z - z_1| \leq \frac{1}{10}(\rh^*)^{1/2}$ in which case 
\begin{equation}
    \begin{split}
        \left|(\rh^*)^{1/2}\sqrt{1+\frac{1}{3}y^2}-\sqrt{1+\frac{1}{3}y_1^2}\right| &= \frac{|\rh^*-1+\frac{1}{3}((\rh^*)^{1/2} y-y_1)((\rh^*)^{1/2} y+ y_1)|}{(\rh^*)^{1/2}\sqrt{1+\frac{1}{3}y^2}+ \sqrt{1+\frac{1}{3}y_1^2}} \\
        &\geq \frac{\rh^*-1-\frac{1}{3}|(\rh^*)^{1/2} y- y_1|(2|(\rh^*)^{1/2} y|+|(\rh^*)^{1/2} y- y_1|)}{\sqrt{4/3}(\rh^*)^{1/2}+\sqrt{1+\frac{1}{3}((\rh^*)^{1/2} y+|(\rh^*)^{1/2} y- y_1|)^2}}\\
        &\geq \frac{\rh^*-1-\frac{\rh^*}{30}(2+1/10)}{(\sqrt{4/3}+\sqrt{1/4+\frac{1}{3}(1+1/10)^2})(\rh^*)^{1/2}} \geq C (\rh^*)^{1/2} \geq C \frac{(\rh^*)^{1/2}}{|z_1|}
    \end{split}
\end{equation}

\end{proof}

Now we're ready for the proof,\\

\textit{Proof of Lemma \ref{Error Lemma for n_1=1}}
Recall the goal is to bound $\mathcal{I}_1[\mathcal{I}_2[\xi_1^{-a}\eta_1^{-b}(\boldsymbol{\xi}^{(2)})^{\alpha}(\boldsymbol{\eta}^{(2)})^{\beta} (\xi_1-\boldsymbol{\xi}^{(2)})^{-\gamma}(\eta_1-\boldsymbol{\eta}^{(2)})^{-\delta}]]$, where $a,b \in \mathbb{N}$ $\alpha,\beta,\gamma,\delta \in \mathbb{N}^{n_2}$ are all multi-indices of length $n_2$, and
\begin{equation}
    \begin{split}
        &\mathcal{I}_2[\xi_1^{-a}\eta_1^{-b}(\boldsymbol{\xi}^{(2)})^{\alpha}(\boldsymbol{\eta}^{(2)})^{\beta} (\xi_1-\boldsymbol{\xi}^{(2)})^{-\gamma}(\eta_1-\boldsymbol{\eta}^{(2)})^{-\delta}]\\
        &:=\prod_{i_2=1}^{n_2}\int_{\Gamma_{\LL,\out}}\ddbar{\xi_{i_2}^{(2)}}{}\int_{\Gamma_{\RR,\out}}\ddbar{\eta_{i_2}^{(2)}}{} \left(C(\boldsymbol{\xi}^{(2)};\boldsymbol{\eta}^{(2)})\right)^2\mathrm{f}_2(\boldsymbol{\xi}^{(2)}+x'/\tau')\mathrm{f}_2(\boldsymbol{\eta}^{(2)}+x'/\tau')\\
        &\xi_1^{-a}\eta_1^{-b}(\boldsymbol{\xi}^{(2)})^{\alpha}(\boldsymbol{\eta}^{(2)})^{\beta} (\xi_1-\boldsymbol{\xi}^{(2)})^{-\gamma}(\eta_1-\boldsymbol{\eta}^{(2)})^{-\delta}
    \end{split}
\end{equation}
Towards this goal we start by making the change of variables $\xi_1=(\rh^*)^{1/2}z$, 
$\eta_1=(\rh^*)^{1/2}z$,
$\xi_{i}^{(2)}=z_i$,
$\eta_{i}^{(2)}=w_{i}$ for $i=1,...,n_2$. And deform all the contours to the steepest descent contours. It's important to remark that no contours need to cross in this change of variables and deformation which would not be true if we were dealing with $\Gamma_{\LL, \inn}$ and $\Gamma_{\RR, \inn}$.\\

Now for $\boldsymbol{z}=(z_1,...,z_{n_2})$, $\boldsymbol{w}=(w_1,...,w_{n_2})$, by lemma \ref{contour difference lemma} on the steepest descent contours,
\begin{equation}
\begin{split}
    &|(\xi_1-\boldsymbol{\xi}_{i_2}^{(2)})^{-\gamma}(\eta_1-\boldsymbol{\eta}_{i_2}^{(2)})^{-\delta}|\\ &=|((\rh^*)^{1/2}z-\boldsymbol{z})^{-\gamma}((\rh^*)^{1/2}w-\boldsymbol{w})^{-\delta}|\\
    & \leq C^{|\gamma|+|\delta|}(\rh^*)^{-\frac{1}{2}(|\gamma|+|\delta|)} |\boldsymbol{z}|^{\gamma}|\boldsymbol{w}|^{\delta}\\
    &\lesssim_{|\gamma|,|\delta|}(\rh^*)^{-\frac{1}{2}(|\gamma|+|\delta|)}|\boldsymbol{z}|^{\gamma}|\boldsymbol{w}|^{\delta}
\end{split}
\end{equation}
The second to last inequality follows from (\ref{distance between contours}) since $|z|,|w|\leq 2$. 

Next, we can use Hadamard's inequality on the Cauchy determinant 
\begin{equation}
    C(\boldsymbol{z};\boldsymbol{w}) = (-1)^{n_2(n_2-1)/2} \mathrm{det}\left[\frac{1}{z_i-w_j}\right]_{i,j=1}^{n_2}
\end{equation} to obtain the bound
\begin{equation}
    \begin{split}
        \left(C(\boldsymbol{\xi}^{(2)};\boldsymbol{\eta}^{(2)})\right)^2 = \left(C(\boldsymbol{z};\boldsymbol{w})\right)^2 \leq n_2^{n_2}
    \end{split}
\end{equation}

So we have 
\begin{equation}
    \begin{split}
        &|\mathcal{I}_1[\mathcal{I}_2[\xi_1^{-a}\eta_1^{-b}(\xi_{i_2}^{(2)})^{\alpha}(\eta_{i_2}^{(2)})^{\beta} (\xi_1-\boldsymbol{\xi}_{i_2}^{(2)})^{-\gamma}(\eta_1-\boldsymbol{\eta}_{i_2}^{(2)})^{-\delta}]]|\\
        &\leq (C'')^{n_2} n_2^{n_2}(\rh^*)^{-\frac{1}{2}(|\gamma|+|\delta|)} J_1 \cdot \prod_{i=1}^{n_2} J_i^{(2)}
    \end{split}
    \label{Formula for I_1[I_2[...]]}
\end{equation}

where 
\begin{equation}
    J_1 = \frac{1}{4\pi^2 (\rh^*)^{(a+b)/2}} \int_{\Gamma_{\LL}^{SD}} \int_{\Gamma_{\RR}^{SD}} \left|\exp(\tau'(\rh^*)^{3/2}((w^3-z^3)/3-(w-z)))\right| \frac{1}{|w-z|^2} |z|^{-a} |w|^{-b} |\dif{w}| |\dif{z}|
\end{equation}
and
\begin{equation}
    J_i^{(2)} = \frac{(\Tilde{h})^{(\alpha_i+\beta_i+2)/2}}{4\pi^2 } \int_{\Gamma_{\LL}^{SD}} \int_{\Gamma_{\RR}^{SD}} \mathrm{f}_2(z_i+x'/\tau')\mathrm{f}_2(w_i+x'/\tau')|z_i|^{\alpha_i+\gamma_i} |w_i|^{\beta_i+\delta_i} |\dif{w_i}| |\dif{z_i}|
\end{equation}

Observe that $|w|,|z|,|w-z| \geq 1$ and $\mathrm{Re}[(w^3-z^3)/3-(w-z)] \leq -4/3$ for $z \in \Gamma_{\LL}^{SD}$, $w \in \Gamma_{\RR}^{SD}$, so for $\tau'(\rh^*)^{3/2}>1$, we have 
\begin{equation}
\begin{split}
    &J_1 \leq \frac{1}{4\pi^2 (\rh^*)^{(a+b)/2}}e^{-\frac{4}{3}\tau'(\rh^*)^{3/2}} \int_{\Gamma_{\LL}^{SD}} \int_{\Gamma_{\RR}^{SD}} \exp(\tau'(\rh^*)^{3/2}(4/3+(w^3-z^3)/3-(w-z))) |\dif{w}| |\dif{z}|\\
    &\leq \frac{1}{4\pi^2 (\rh^*)^{(a+b)/2}}e^{-\frac{4}{3}\tau'(\rh^*)^{3/2}} \int_{\Gamma_{\LL}^{SD}} \int_{\Gamma_{\RR}^{SD}} \exp((4/3+(w^3-z^3)/3-(w-z))) |\dif{w}| |\dif{z}|\\
    &\leq C_1 (\rh^*)^{-(a+b)/2} e^{-\frac{4}{3}\tau'(\rh^*)^{3/2}}
    \end{split}
    \label{I_1 bound}
\end{equation}

In order to estimate $J_i^{(2)}$, we follow an identical argument to \eqref{I_2 of monomials bound}. Observe that $\mathrm{f}_2(z_i+x'/\tau')\mathrm{f}_2(w_i+x'/\tau') \lesssim e^{-\frac{1}{4}\tau(|z_i|^3+|w_i|^3)}$ for $w_i,z_i$ on the steepest descent contours, thus

\begin{equation}
    \begin{split}
        J_i^{(2)} &\lesssim\frac{1}{4\pi^2 }  \int_{\Gamma_{\LL}^{SD}} \int_{\Gamma_{\RR}^{SD}} e^{-\frac{1}{4}\tau (|z_i|^3+|w_i|^3)}|z_i|^{\alpha_i+\gamma_i} |w_i|^{\beta_i+\delta_i} |\dif{w_i}| |\dif{z_i}|\\&\lesssim\frac{1}{4\pi^2 }  \int_{\Gamma_{\LL}^{SD}} \int_{\Gamma_{\RR}^{SD}} e^{-\frac{1}{4}\tau (|z_i|^3+|w_i|^3)}|z_i|^{|\alpha|_{\infty}+|\gamma|_{\infty}} |w_i|^{|\beta|_{\infty}+|\delta|_{\infty}} |\dif{w_i}| |\dif{z_i}|\\
        &\lesssim_{|\alpha|_{\infty},|\beta|_{\infty},|\gamma|_{\infty},|\delta|_{\infty}} 1
    \end{split}
    \label{I_i^(2) bound}
\end{equation}

Finally combining \eqref{Formula for I_1[I_2[...]]}, \eqref{I_1 bound}, and \eqref{I_i^(2) bound}, we get the desired bound.

\end{document}